\nonstopmode
\documentclass[a4paper, 10pt, reqno]{amsart}
\usepackage{latexsym}
\usepackage{fancyhdr}
\usepackage{amsmath, amssymb}
\usepackage[ansinew]{inputenc}
\usepackage[all]{xy}
\theoremstyle{plain}
\newtheorem*{lemma*}{Lemma}
\newtheorem{lemma}[subsection]{Lemma}
\newtheorem*{theorem*}{Theorem}
\newtheorem{theorem}[subsection]{Theorem}
\newtheorem*{proposition*}{Proposition}
\newtheorem{proposition}[subsection]{Proposition}
\newtheorem*{corollary*}{Corollary}
\newtheorem{corollary}[subsection]{Corollary}
\theoremstyle{definition}
\newtheorem*{definition*}{Definition}

\newtheorem*{example*}{Example}
\newtheorem{example}[subsection]{Example}
\theoremstyle{remark}
\newtheorem*{remark*}{Remark}
\newtheorem{remark}[subsection]{Remark}

\numberwithin{equation}{section}

\def\B{\mathcal{B}}
\def\id{\mathrm{id}}

\def\pr{\mathrm{pr}}
\def\reg{\mathrm{reg}}

\title[Lifting smooth curves over invariants, III]
{Lifting smooth curves over invariants\\
for representations of compact Lie groups, III}

\author[A. Kriegl, M. Losik, P.W. Michor, A. Rainer]
{Andreas Kriegl, Mark Losik, Peter W. Michor, and Armin Rainer}

\address{A. 
Kriegl: Institut f\"ur Mathematik, Universit\"at Wien, 
Nordbergstrasse~15, A-1090 Wien, Austria}

\email{Andreas.Kriegl@univie.ac.at}

\address{M. 
Losik: Saratov State University, 
ul. 
Astrakhanskaya, 83, 410026 Saratov, Russia}

\email{losikMV@info.sgu.ru}

\address{P.W. 
Michor: Institut f\"ur Mathematik, Universit\"at Wien, 
Nordbergstrasse~15, A-1090 Wien, Austria; and:
Erwin Schr\"odinger Institut f\"ur Mathematische Physik, 
Boltzmann\-gasse 9, A-1090 Wien, Austria}

\email{Peter.Michor@esi.ac.at}

\address{A. 
Rainer: Institut f\"ur Mathematik, Universit\"at Wien, 
Nordbergstrasse~15, A-1090 Wien, Austria}

\email{armin.rainer@univie.ac.at}

\begin{document}

\begin{abstract}
Any sufficiently often differentiable curve in the orbit space $V/G$ 
of a real finite-dimensional orthogonal representation $G \to O(V)$ of a 
finite group $G$ admits a differentiable lift into the representation space $V$ 
with locally bounded derivative.
As a consequence any sufficiently often differentiable curve in the orbit 
space $V/G$ can be lifted twice differentiably.  
These results can be generalized to arbitrary polar representations.
Finite reflection groups and finite rotation groups in dimensions two and 
three are discussed in detail.
\end{abstract}

%\begin{center}
%Erwin Schr\"odinger International Institute \\
%of Mathematical Physics, Wien, Austria
%\end{center}
\thanks{M.L., P.W.M. and A.R. were supported by
`Fonds zur F\"orderung der wissenschaftlichen Forschung, Projekt P 17108-N04'}
%\dedicatory{}
\date {\today}

\maketitle

\section{Introduction}

In \cite{rep-lift} and \cite{rep-lift3} the following problem was investigated.
Consider an orthogonal representation of a compact Lie group $G$
on a real finite dimensional Euclidean vector space $V$.
Let $\sigma_1,\ldots,\sigma_n$ be a system of homogeneous generators
for the algebra $\mathbb{R}[V]^G$ of invariant polynomials on $V$.
Then the mapping
$\sigma = (\sigma_1,\ldots,\sigma_n) : V \rightarrow \mathbb{R}^n$
induces a homeomorphism between the orbit space $V/G$
and the semialgebraic set $\sigma(V)$.
Suppose a smooth curve
$c : \mathbb{R} \rightarrow V/G = \sigma(V) \subseteq \mathbb{R}^n$
in the orbit space is given (smooth as curve in $\mathbb{R}^n$),
does there exist a smooth lift to $V$,
i.e., a smooth curve $\bar{c} : \mathbb{R} \rightarrow V$ with
$c = \sigma \circ \bar{c}$\,? The answer is independent of the choice of
the generators.

It was shown in \cite{rep-lift} that a real analytic curve in $V/G$
admits a local real analytic lift to $V$, and that
a smooth curve in $V/G$ admits a global smooth lift,
if certain genericity conditions are satisfied.
In both cases the lifts
may be chosen orthogonal to each orbit they meet and then they are
unique up to a transformation in $G$,
whenever the representation of $G$ on $V$ is polar,
i.e., admits sections.

In \cite{rep-lift3} we proved that any sufficiently often differentiable 
curve in the orbit space $V/G$ can be lifted to a once differentiable 
curve in $V$.

In the special case that the symmetric group $S_n$ is acting on 
$\mathbb R^n$ by permuting the coordinates there is the following 
interpretation of the described lifting problem. 
As generators of $\mathbb{R}[\mathbb R^n]^{S_n}$ we may take the 
elementary symmetric functions
\[
\sigma_j(x) = 
\sum_{1 \le i_1 < \cdots < i_j \le n} x_{i_1} \cdots x_{i_j}
\qquad (1 \le j \le n),
\]
which constitute the coefficients (up to sign) of a monic polynomial $P$ 
with roots $x_1,\ldots,x_n$ via Vieta's formulas.
Then a curve in the orbit space $\mathbb R^n/S_n = \sigma(\mathbb R^n)$
corresponds to a curve $P(t)$ of monic polynomials of degree $n$ with 
only real roots, and a lift of $P(t)$ may be interpreted as a 
parameterization of the roots of $P(t)$. 

This problem has been studied extensively in \cite{roots}.
Moreover, the following results were proved in \cite{rootsII}:
Any differentiable lift (parameterization of roots) of a $C^{2 n}$-curve 
(of polynomials) $P : \mathbb{R} \rightarrow \mathbb{R}^n/S_n$ is 
actually $C^1$,
and there always exists a twice differentiable
but in general not better lift of $P$,
if $P$ is of class $C^{3 n}$.
Note that here the differentiability assumptions on $P$ are not the
weakest possible which is shown by the case $n=2$,
elaborated in \cite{roots} 2.1.
The proof in \cite{rootsII} is based on the fact that the roots of a
$C^n$-curve of
polynomials $P : \mathbb{R} \rightarrow \mathbb{R}^n/S_n$
may be chosen differentiable with locally bounded derivative;
this is due to Bronshtein \cite{bronshtein} and
Wakabayashi \cite{wakabayashi}.

In the present paper we show the corresponding statements for arbitrary 
real finite-dimensional orthogonal representations of finite groups. 
We consider representations $\rho : G \to O(V)$ with 
the \emph{property $(\B_k)$} that any $C^k$-curve in a neighborhood of $0$ in the 
orbit space $V/G$ admits a local differentiable lift to $V$ 
with locally bounded derivative (section \ref{secB}). 
In analogy to the polynomial case 
$S_n : \mathbb{R}^n$ we then show that, for representations of finite groups 
$G$ with property $(\B_k)$, any differentiable lift of a $C^{k+d}$-curve 
is actually $C^1$ (section \ref{secc1}), and there always exists a 
twice differentiable lift of a $C^{k+2d}$-curve in the orbit space 
(section \ref{sec2der}). 
The integer $d$ is the maximal degree of a minimal system of homogeneous 
generators of the algebra of invariant polynomials $\mathbb R[V]^G$ 
(see \ref{dD}).
As a consequence we obtain in section \ref{secpolar} 
that polar representations, where the representation of the 
associated generalized Weyl group on some section has property $(\B_k)$, 
allow orthogonal $C^1$-lifts of $C^{k+d}$-curves and orthogonal 
twice differentiable lifts of $C^{k+2d}$-curves. 
In section \ref{secstab} we show that property $(\B_k)$ is stable 
under subrepresentations and orthogonal direct sums. 
We prove in section \ref{secallfiniteG}, by reducing to the polynomial case, 
that any real finite-dimensional 
representation $\rho : G \to O(V)$ of a finite group $G$ has property $(\B_k)$, 
where 
\[
k = \max\{d, |G|/|G_{v_i}| : 1 \le i \le l\},
\] 
$v_i \in V_i \backslash \{0\}$ are chosen such that the cardinality of the 
isotropy groups $G_{v_i}$ is maximal, 
and $V = V_1 \oplus \cdots \oplus V_l$ is the 
decomposition into irreducible subrepresentations. 
This establishes property $(\B)$ for polar representations, too. 
In section \ref{reflectiongroups} we give a complete survey of all finite 
reflection groups. 
Section \ref{rotgroups} is devoted to the discussion of finite rotation 
groups in dimensions two and three.

Still open is the question whether non-polar representations 
of compact Lie groups $G$ on real finite dimensional Euclidean 
vector spaces $V$ have property $(\B)$.

The polynomial results have applications in the theory
of partial differential equations and perturbation
theory, see \cite{perturb}.

\section{Preliminaries}

\subsection{The setting} \label{setting}

Let $G$ be a compact Lie group and let $\rho : G \rightarrow O(V)$ be an
orthogonal representation in a real finite dimensional Euclidean
vector space $V$ with inner product $\langle \quad |\quad \rangle$.
By a classical theorem of Hilbert and Nagata,
the algebra $\mathbb{R}[V]^{G}$ of invariant polynomials on $V$
is finitely generated.
So let $\sigma_1,\ldots,\sigma_n$ be a system of homogeneous generators
of $\mathbb{R}[V]^{G}$ of positive degrees $d_1,\ldots,d_n$.
We may assume that $\sigma_1 : v \mapsto \langle v|v \rangle$
is the inner product.
Consider the \emph{orbit map}
$\sigma = (\sigma_1,\ldots,\sigma_n) : V \rightarrow \mathbb{R}^n$.
Note that,
if $(y_1,\ldots,y_n) = \sigma(v)$ for $v \in V$, then
$(t^{d_1} y_1,\ldots,t^{d_n} y_n) = \sigma(t v)$ for $t \in \mathbb{R}$,
and that $\sigma^{-1}(0) = \{0\}$.
The image $\sigma(V)$ is a semialgebraic set in the categorical quotient
$V/\!\!/G := \{y \in \mathbb{R}^n : P(y) = 0 ~\mbox{for all}~ P \in I\}$
where $I$ is the ideal of relations between $\sigma_1,\ldots,\sigma_n$.
Since $G$ is compact, $\sigma$ is proper and separates orbits of $G$,
it thus induces a homeomorphism between $V/G$ and $\sigma(V)$.

\subsection{The problem of lifting curves} \label{liftingproblem}
Let $c : \mathbb{R} \to V/G = \sigma(V) \subseteq \mathbb{R}^n$ be a smooth curve 
in the orbit space; smooth as curve in $\mathbb{R}^n$.
A curve $\bar{c} : \mathbb{R} \to V$ is called lift of $c$ to $V$, if $c = \sigma \circ \bar{c}$ 
holds. 
\emph{The problem of lifting smooth curves over invariants is independent of the 
choice of a system of homogeneous generators of $\mathbb{R}[V]^G$} 
in the following sense: Suppose $\sigma_1,\ldots,\sigma_n$ and 
$\tau_1,\ldots,\tau_m$ both generate $\mathbb{R}[V]^G$. 
Then for all $i$ and $j$ we have 
$\sigma_i = p_i(\tau_1,\ldots,\tau_m)$ and 
$\tau_j = q_j(\sigma_1,\ldots,\sigma_n)$ for polynomials $p_i$ and $q_j$. 
If $c^{\sigma} = (c_1,\ldots,c_n)$ is a curve in $\sigma(V)$, then 
$c^{\tau} = (q_1(c^{\sigma}),\ldots,q_m(c^\sigma))$ defines a curve in 
$\tau(V)$ of the same regularity. Any lift $\bar{c}$ to $V$ of the curve 
$c^{\sigma}$, i.e., $c^{\sigma} = \sigma \circ \bar{c}$, is a lift of 
$c^{\tau}$ as well (and conversely):
\[c^{\tau} = (q_1(c^{\sigma}),\ldots,q_m(c^\sigma))
=(q_1(\sigma(\bar{c})),\ldots,q_m(\sigma(\bar{c})))
=(\tau_1(\bar{c}),\ldots,\tau_m(\bar{c})) 
=\tau \circ \bar{c}.
\]
%Obviously, the same is true for the converse.

\subsection{The slice theorem} \label{slice}
For a point $v \in V$ we denote by $G_v$ its isotropy group and by
$N_v = T_v(G.v)^{\bot}$ the normal subspace of the orbit $G.v$ at $v$.
It is well known that there exists a $G$-invariant open neighborhood $U$
of $v$ which is real analytically $G$-isomorphic to the crossed product
(or associated bundle) $G \times_{G_v} S_v = (G \times S_v)/G_v$,
where $S_v$ is a ball in $N_v$ with center at the origin.
The quotient $U/G$ is homeomorphic to $S_v/G_v$.
It follows that the problem of local lifting curves in $V/G$ passing
through $\sigma(v)$ reduces to the same problem for curves in $N_v/G_v$
passing through $0$.
For more details see \cite{rep-lift}, \cite{luna}, and \cite{schwarz1}
theorem 1.1.

A point $v \in V$ (and its orbit $G.v$ in $V/G$) is called \emph{regular}
if the slice representation
$G_v \rightarrow O(N_v)$ is trivial.
Hence a neighborhood of this point is analytically $G$-isomorphic to
$G/G_v \times S_v \cong G.v \times S_v$.
The set $V_{\reg}$ of regular points is open and dense in $V$,
and the projection $V_{\reg} \rightarrow V_{\reg}/G$
is a locally trivial fiber bundle.
A non regular orbit or point is called \emph{singular}.

\subsection{The integer $d$} \label{dD}
Let $\rho : G \to O(V)$ be as in \ref{setting}. 
Choose a minimal system of homogeneous generators $\sigma_1,\ldots,\sigma_n$ 
of positive degrees $d_1,\ldots,d_n$ of $\mathbb{R}[V]^G$. 
We associate to $\rho$ the following number: 
\[d = d(\rho) := \max \{d_1,\ldots,d_n\}.\]
\emph{The integer $d$ is well-defined and independent of the choice of a 
minimal system of homogeneous generators of the algebra of invariant polynomials.} 
This follows from the fact that a system of homogeneous invariants 
of positive degree generates $\mathbb{R}[V]^G$ as an algebra over $\mathbb{R}$ 
if and only if the images of the invariants in this system generate 
$\mathbb{R}[V]^G_+/(\mathbb{R}[V]^G_+)^2$ as a vector space over $\mathbb{R}$, 
where $\mathbb{R}[V]^G_+$ is the space of all invariants with positive degree, 
e.g. \cite{DerksenKemper}, 3.6. 
The grading used here is given by the degree of the polynomials. 
Hence a system of homogeneous algebra generators has minimal cardinality 
if no generator is superfluous, 
and the number and the degrees of the elements in a minimal system of 
homogeneous generators are uniquely determined.

Note that independence of $d$ from the choice of a minimal system of 
homogeneous generators of $\mathbb{R}[V]^G$ also follows from the following lemma 
applied to the slice representation at $0$.

\begin{lemma*} 
Let $\rho : G \to O(V)$ be a finite-dimensional representation of a 
compact Lie group $G$, 
let $\rho'$ be some slice representation of $\rho$. 
Then, $d(\rho') \le d(\rho)$.  
\end{lemma*}

\proof
Let $\sigma_1,\ldots,\sigma_n$ be a minimal system of homogeneous generators 
of $\mathbb{R}[V]^G$.

For an arbitrary $v \in V$ let $\rho' : G_v \to O(N_v)$ be its slice 
representation, and suppose $S_v$ is a normal slice at $v$. 
Choose a minimal system of homogeneous generators $\tau_1,\ldots,\tau_m$ of 
$\mathbb{R}[N_v]^{G_v}$ and assume that 
$\deg \tau_1 \le \cdots \le \deg \tau_m = d(\rho')$. 
Then there exist polynomials $p_i \in \mathbb{R}[\mathbb{R}^m]$ such that 
\[
\sigma_i|_{S_v} = p_i(\tau_1|_{S_v},\ldots,\tau_m|_{S_v}) \qquad (1 \le i \le n).
\]
On the other hand, by the slice theorem, near $v \in N_v$ we have 
\[
\tau_j|_{S_v} = f_j(\sigma_1|_{S_v},\ldots,\sigma_n|_{S_v}) \qquad (1 \le j \le m),
\]
where $f_j$ are real analytic functions; e.g. \cite{schwarz2}.

For contradiction assume $\deg \tau_m > d(\rho)$. 
Then all polynomials $p_i$ do not depend on their last entry. 
Consequently, near $v \in N_v$,
\[
\tau_m|_{S_v} = F(\tau_1|_{S_v},\ldots,\tau_{m-1}|_{S_v}),
\]
where
\[
F = f_m(p_1,\ldots,p_n)
\]
is real analytic. 
Introduce a new grading in $\mathbb{R}[\mathbb{R}^{m-1}]$ with respect to 
$\deg \tau_1 \le \cdots \le \deg \tau_m$ and write the function $F$ as an 
infinite sum of homogeneous (with respect to this grading) terms. 
Let $\bar{F}$ be the sum of all terms of degree $\deg \tau_m$ in this 
presentation of $F$. 
We obtain, near $v \in N_v$, 
\[
\tau_m|_{S_v} = \bar{F}(\tau_1|_{S_v},\ldots,\tau_{m-1}|_{S_v}).
\]
This means $\tau_m$ is a polynomial in $\tau_1,\ldots,\tau_{m-1}$ in a 
neighborhood of $v$ in $N_v$, and, hence, everywhere. 
This contradicts minimality of $\tau_1,\ldots,\tau_m$.
\endproof

\begin{remark*}
The previous lemma allows to replace the intricate definition of the integer $d$ 
given in \cite{rep-lift3} by the definition given above. 
\end{remark*}

\subsection{Removing fixed points} \label{fix}

Let $V^G$ be the space of $G$-invariant vectors in $V$,
and let $V'$ be its orthogonal complement in $V$.
Then we have $V = V^G \oplus V'$,
$\mathbb{R}[V]^G = \mathbb{R}[V^G] \otimes \mathbb{R}[V']^G$, and
$V/G = V^G \times V'/G$.

\begin{lemma*}
Any lift $\bar{c}$ of a curve $c=(c_0,c_1)$ of class $C^k$
($k = 0,1,\ldots,\infty,\omega$) in $V^G \times V'/G$ has the form
$\bar{c} = (c_0,\bar{c}_1)$,
where $\bar{c}_1$ is a lift of $c_1$ to $V'$ of class
$C^k$.
Then the lift $\bar{c}$ is orthogonal if and only if the lift
$\bar{c}_1$ is orthogonal. \qed
\end{lemma*}

\subsection{Multiplicity} \label{mult}

For a continuous function $f$ defined near $0$ in $\mathbb{R}$,
let the \emph{multiplicity} or \emph{order of flatness} $m(f)$ at
$0$ be the supremum of all integers $p$ such that $f(t) = t^p g(t)$ near
$0$ for a continuous function $g$.
If $f$ is $C^n$ and $m(f) < n$,
then $f(t) = t^{m(f)} g(t)$,
where now $g$ is $C^{n-m(f)}$ and $g(0) \ne 0$.
Similarly,
one can define multiplicity of a function at any $t \in \mathbb{R}$.

\begin{lemma*}
Let $c = (c_1,\ldots,c_n)$ be a curve in $\sigma(V) \subseteq \mathbb{R}^n$ 
of class $C^r$, where $r \ge d$,
and $c(0) = 0$.
Then the following two conditions are equivalent:
\begin{enumerate}
\item $c_1(t) = t^2 c_{1,1}(t)$ near $0$ for a $C^{r-2}$-function $c_{1,1}$;
\item $c_i(t) = t^{d_i} c_{i,i}(t)$ near $0$ 
for a $C^{r-d_i}$-function $c_{i,i}$, for all $1 \le i \le n$.
\end{enumerate}
\end{lemma*}

\proof
The proof is the same as that of lemma 2.4 in \cite{rep-lift3}.\qed

\subsection{} \label{recall}

We recall a few facts from \cite{rep-lift3}:

\begin{lemma*}[a]
A curve
$c : \mathbb{R} \rightarrow V/G = \sigma(V) \subseteq \mathbb{R}^n$
of class $C^d$ admits an orthogonal $C^d$-lift $\bar{c}$ in a
neighborhood of a regular point $c(t_0) \in V_{\reg}/G$.
It is unique up to a transformation from $G$.
\end{lemma*}

\begin{lemma*}[b]
Consider a continuous curve $c : (a,b) \to X$ in a compact metric space $X$. 
Then the set $A$ of all accumulation points of $c(t)$ as 
$t \searrow a$ is connected.
\end{lemma*}

\begin{theorem*}
Let
$c = (c_1,\ldots,c_n) : \mathbb{R} \rightarrow V/G
  = \sigma(V) \subseteq \mathbb{R}^n$
be a curve of class $C^d$.
Then there exists a global differentiable lift
$\bar{c} : \mathbb{R} \rightarrow V$ of $c$.
\end{theorem*}

\subsection{Bronshtein's and Wakabayashi's result} \label{bron-wa}

We formulate Bronshtein's theorem \cite{bronshtein}; 
consult also Wakabayashi's version \cite{wakabayashi}.

\begin{theorem*} \cite{bronshtein}
Let
\[P(t)(x) = x^n + \sum_{j=1}^n (-1)^j a_j(t) x^{n-j}\]
be a curve of monic polynomials of degree $n$ with all roots real for all 
$t \in \mathbb{R}$, where $a_j \in C^n(\mathbb{R})$ for all $1 \le j \le n$. 
Choose a differentiable parameterization $x_1(t),\ldots,x_n(t)$ 
of the roots of $P(t)$ (which always exists). Then, for any compact subset 
$K \subseteq \mathbb{R}$ 
there exists a constant $C_K$ such that 
\[\left|\frac{d}{d t} x_j(t)\right| \le C_K \qquad \mbox{for all} ~ t \in K, 
1 \le j \le n.\]

In the language of representation theory: Any $C^n$-curve $P$ in the orbit 
space $\mathbb{R}^n/S_n$ of the standard representation of the 
symmetric group $S_n$ on $\mathbb{R}^n$, by permuting the coordinates, 
allows a differentiable lift $x=(x_1,\ldots,x_n)$ 
with locally bounded derivative. 
\end{theorem*}

\section{Property $(\B)$} \label{secB}

\subsection{Property $(\B)$} \label{propertyB}
We shall say that an orthogonal representation $\rho : G \to O(V)$ 
of a compact Lie group $G$ on a real finite dimensional 
Euclidean vector space $V$ 
has \emph{property $(\B_k)$}, if:

\begin{quote}
There exists a neighborhood $U = U(\rho)$ of $0$ in $V/G = \sigma(V)$ 
such that each $C^k$-curve 
in $U$ admits a local differentiable lift $\bar{c}$ to $V$ with locally 
bounded derivative.
\end{quote}

Note that property $(\B_k)$ is independent of the choice of generators of 
$\mathbb{R}[V]^G$.

It is clear that, if a representation $\rho$ has property $(\B_k)$, 
then it has property $(\B_l)$ for all $l \in \{k,k+1,\ldots,\infty,\omega\}$ 
as well.

We shall write simply property $(\B)$, if the degree of differentiability $k$ 
is not specified.

\begin{example*}
The standard representation of the symmetric group $S_n$ on $\mathbb{R}^n$ 
has property $(\B_n)$. This follows from theorem \ref{bron-wa}.
\end{example*}

\begin{proposition} \label{Beverywhere}
Let
$c = (c_1,\ldots,c_n) : \mathbb{R} \rightarrow V/G
= \sigma(V) \subseteq \mathbb{R}^n$ be a curve of class $C^{k}$ in the 
orbit space of a representation $\rho : G \to O(V)$ with property $(\B_k)$.
Then for any $t_0 \in \mathbb{R}$ there exists a local differentiable lift
$\bar{c}$ of $c$ near $t_0$ with locally bounded derivative.
\end{proposition}

\proof
For each $s \in \mathbb{R} \backslash \{0\}$ let us define a $C^k$-curve 
$c_s : \mathbb{R} \to \sigma(V)$ by 
\[
c_s(t) = (s^{d_1} c_1(t),\ldots,s^{d_n} c_n(t)).
\]
There exists some $s=s(c;t_0) \in \mathbb{R} \backslash \{0\}$ such that 
$c_s(t) \in U$ for $t$ near $t_0$, where $U$ is the neighborhood of $0$ in 
$V/G$ introduced in the definition of property $(\B_k)$. 
Since $\rho$ has property $(\B_k)$, there exists, near $t_0$, a local 
differentiable lift $\bar{c}_s$ of $c_s$ to $V$ with locally bounded derivative. 
Then, $\bar{c}(t) := s^{-1} \cdot \bar{c}_s(t)$ defines a local differentiable 
lift of $c$ for $t$ near $t_0$ whose derivative is locally bounded.
\endproof

\begin{proposition} \label{Bk&slice}
Assume that $\rho : G \to O(V)$ is a representation of a 
finite group $G$ with property $(\B_k)$. Then any slice representation 
$\rho'$ of $\rho$ has property $(\B_k)$ as well.
\end{proposition}

\proof
Let $\rho' : G_v \to O(N_v)$ be an arbitrary slice representation of $\rho$. 
Consider some normal slice $S_v$ at $v$ for the $G$-action on $V$. 
Then $S_v/G_v$ is an open neighborhood of $0$ in $N_v/G_v$ which by \ref{slice} 
is homeomorphic to $(G \times_{G_v} S_v)/G$ which in turn is an open neighborhood 
of $G.v$ in $V/G$. 

Given a $C^k$-curve $c$ in $S_v/G_v$, we may view it as a curve in 
$(G \times_{G_v} S_v)/G$. Since $\rho$ has property $(\B_k)$ and by proposition 
\ref{Beverywhere}, there exists a local differentiable lift $\bar{c}$ of $c$ to 
$V$ with locally bounded derivative. 
The finiteness of $G$ implies that $N_v = V$, and hence $S_v$ is an open 
neighborhood of $v$ in $V$. 
Therefore $\bar{c}$ is a local lift of $c$ to $N_v$ with respect to the $G_v$-action.
\endproof

\begin{lemma} \label{deroflifts}
Let
$c : \mathbb{R} \rightarrow V/G
= \sigma(V) \subseteq \mathbb{R}^n$ be a curve in the orbit space $V/G$. 
We assume that $G$ is finite. Let $t_0 \in \mathbb{R}$. 
If $\bar{c}_1$ and $\bar{c}_2$ are lifts of $c$ which are (one-sided) 
differentiable at $t_0$ and $\bar{c}_1(t_0) = \bar{c}_2(t_0)$, 
then there exists some $g \in G_{\bar{c}_1(t_0)}$ 
such that $\bar{c}_1'(t_0) = g.\bar{c}_2'(t_0)$.
\end{lemma}

\proof
Without loss we can assume that $t_0 = 0$. 

Let $\bar{c}_1$ and $\bar{c}_2$ be lifts of $c : \mathbb{R} \to V/G$ 
which are (one-sided) differentiable at $0$ and satisfy 
$\bar{c}_1(0) = \bar{c}_2(0) =: v_0$. 
We may suppose $V^G = \{0\}$, by lemma \ref{fix}. 
We consider the following cases separately:

If $c(0) = 0$, then $\bar{c}_1(0) = \bar{c}_2(0) = 0$ and consequently, 
for $i=1,2$,
\[\sigma(\bar{c}_i'(0)) 
  = \sigma \left( \lim_{t \to 0} \frac{\bar{c}_i(t)}{t}\right) 
  = \lim_{t \to 0} \sigma \left(  \frac{\bar{c}_i(t)}{t}\right) .\]
Now, for $t \ne 0$ we have $\sigma \left( \bar{c}_i(t)/t\right) = c_{(1)}(t) \in \sigma(V)$, 
where 
\[
c_{(1)}(t) := (t^{-d_1} c_1(t),\ldots,t^{-d_n} c_n(t)). 
\]
Since $\sigma(V)$ is closed in $\mathbb{R}^n$
(see \cite{procesischwarz}), 
we find
\[
\sigma(\bar{c}_i'(0)) 
= \lim_{t \to 0} \sigma \left(  \frac{\bar{c}_i(t)}{t}\right)  
= \lim_{t \to 0} c_{(1)}(t) \in \sigma(V),
\]
i.e., $\sigma$ maps $\bar{c}_1'(0)$ and $\bar{c}_2'(0)$ to the same point in $\sigma(V)$. 
(Note that, if only one-sided derivatives exist, then $t \to 0$ has 
to be replaced by $t \nearrow 0$ or $t \searrow 0$, respectively.) 
This shows that $\bar{c}_1'(0)$ and $\bar{c}_2'(0)$ lie in the same orbit, 
therefore we find some $g \in G = G_0$ with $\bar{c}_1'(0) = g.\bar{c}_2'(0)$.

If $c(0) \ne 0$: 
Since $G$ is finite and therefore $N_{v_0} = V$, 
the ball $S_{v_0}$ is a neighborhood of $v_0$ in $V$ 
which contains the lifts $\bar{c}_1(t)$ and $\bar{c}_2(t)$ for $t$ near $0$. 
Hence, by \ref{slice}, we may change to the slice 
representation $G_{v_0} \to O(N_{v_0})$. 
Now we may assume that $c$ is a curve in $N_{v_0}/G_{v_0}$ with $c(0) = 0$ 
and with lifts $\bar{c}_1(t)$ and $\bar{c}_2(t)$ to $N_{v_0}$ for $t$ near $0$. 
So we refer to the former case.
\endproof

Note that lemma \ref{deroflifts} does no longer hold, if finiteness of $G$ 
is omitted: 
\begin{example*}
Consider the standard action of $SO(2)$ on $\mathbb{R}^2$. 
Then $\sigma(x_1,x_2) = x_1^2+x_2^2$ generates 
$\mathbb{R}[\mathbb{R}^2]^{SO(2)}$ and 
$\mathbb{R}^2/SO(2) = \sigma(\mathbb{R}^2) = [0,\infty)$. 
We consider the curve $c(t)=t^2$ and its differentiable lifts 
$\bar{c}_1(t) = (t,0)$ and $\bar{c}_2(t) = (t \cos t,t \sin t)$. 
We find $\bar{c}_1(2 \pi) = \bar{c}_2(2 \pi) =(2 \pi,0)$, 
but $\bar{c}_1'(2 \pi) = (1,0)$ and $\bar{c}_2'(2 \pi) = (1,2 \pi)$ 
cannot be transformed to each other by an element of $G_{(2 \pi,0)}=\{\id\}$.
\end{example*}

\begin{remark*}
If $G$ is not finite, then lemma \ref{deroflifts} generalizes to the 
following statement:
\emph{Let
$c : \mathbb{R} \rightarrow V/G
= \sigma(V) \subseteq \mathbb{R}^n$ be a curve in the orbit space $V/G$. 
Let $t_0 \in \mathbb{R}$. 
If $\bar{c}_1$ and $\bar{c}_2$ are lifts of $c$ which are (one-sided) 
differentiable at $t_0$ and $\bar{c}_1(t_0) = \bar{c}_2(t_0) =:v_0$, 
then there exists some $g \in G_{v_0}$ such that 
$\bar{c}_1'(t_0)^{\bot} = g.\bar{c}_2'(t_0)^{\bot}$, 
where $\bot$ indicates the projection onto $N_{v_0}$.
}

To see this:
We consider the projection 
$p : G.S_{v_0} \cong G \times_{G_{v_0}} S_{v_0} \to G/G_{v_0} \cong G.v_0$ 
of a fiber bundle associated to the principal bundle $\pi : G \to G/G_{v_0}$, 
where $S_{v_0}$ is a normal slice at $v_0$.
Then,
for $t$ close to $t_0$, $\bar{c}_1$ and $\bar{c}_2$ are curves in $G.S_{v_0}$,
whence $p \circ \bar{c}_i$ $(i =1,2)$ are curves in $G/G_{v_0}$ 
which admit lifts $g_i$ into $G$ with $g_i(t_0) = e$, 
which are (one-sided) differentiable at $t_0$ 
(via the horizontal lift of the principal connection, say). 
Consequently, $t \mapsto g_i(t)^{-1}.\bar{c}_i(t)$ are lifts which lie in 
$S_{v_0}$, whence 
$\left.\frac{d}{d t}\right|_{t = t_0} (g_i(t)^{-1}.\bar{c}_i(t)) 
= - g_i'(t_0).v_0 + \bar{c}_i'(t_0) \in N_{v_0}$. 
Thus, 
$\bar{c}_i'(t_0)^{\bot} 
= \left.\frac{d}{d t}\right|_{t = t_0} (g_i(t)^{-1}.\bar{c}_i(t))$. 
By this observation, we may assume without loss that the 
lifts $\bar{c}_1$ and $\bar{c}_2$ lie in $S_{v_0}$ for $t$ close to $t_0$. 
Then the proof of lemma \ref{deroflifts} 
gives the statement.
\end{remark*}

\begin{remark} \label{rmk}
Lemma \ref{deroflifts} implies that for any two differentiable lifts 
$\bar{c}_1$ and $\bar{c}_2$ of a curve $c$ in $V/G$, 
where $G$ is finite, we have $\| \bar{c}_1'(t)\| = \|\bar{c}_2'(t)\|$ 
for all $t$. 
So, if there exists some differentiable lift of $c$ with 
locally bounded derivative, then any differentiable lift of $c$ 
has this property as well.
\end{remark}

\begin{proposition} \label{Bglobal}
Assume that $\rho : G \to O(V)$ is a representation of a finite group $G$ 
with property $(\B_k)$. 
Let $c : \mathbb{R} \to V/G = \sigma(V) \subseteq \mathbb{R}^n$ be a 
curve of class $C^k$.
Then there exists a global differentiable lift $\bar{c}$ of $c$ to $V$ with 
locally bounded derivative.
\end{proposition}

\proof
Proposition \ref{Beverywhere} provides local differentiable lifts of $c$ with 
locally bounded derivative near any $t \in \mathbb{R}$.

Now let us construct from these data a global differentiable lift of $c$ with 
locally bounded derivative:
First we glue the local differentiable lifts with locally bounded derivative 
just differentiably. 
It is sufficient to show that each local differentiable lift of $c$ defined on 
an open interval $I$ can be extended to a larger interval whenever 
$I \ne \mathbb{R}$. 

Suppose that $\bar{c}_1 : I \to V$ is a local differentiable lift of $c$, and 
suppose the open interval $I$ is bounded from above (say), and $t_1$ is its 
upper boundary point. 
Then, there exists a local differentiable lift $\bar{c}_2$ of $c$ 
near $t_1$, and a $t_2 < t_1$ such that both $\bar{c}_1$ and $\bar{c}_2$ are 
defined near $t_2$. 
There is some $g \in G$ such that $\bar{c}_1(t_2) = g.\bar{c}_2(t_2)$. 
By lemma \ref{deroflifts}, we find an $h \in G_{\bar{c}_1(t_2)}$ with 
$\bar{c}_1'(t_2) = hg.\bar{c}_2'(t_2)$. 
Then $\bar{c}(t) := \bar{c}_1(t)$ for $t \le t_2$ and 
$\bar{c}(t) := hg.\bar{c}_2(t)$ for $t \ge t_2$ defines a differentiable 
lift of $c$ on a larger interval. 

Now let us show that the resulting global differentiable lift $\bar{c}$ of 
$c$ has locally bounded derivative. 
Note first that by remark \ref{rmk} the gluing process described above 
does not affect the local boundedness of the derivatives of the local lifts, 
provided by proposition \ref{Beverywhere}. 
Let $K$ be a compact subset of $\mathbb{R}$. 
The domains of definition of the local lifts constitute 
an open covering of $K$ which contains a finite open subcovering 
$\{I_j\}$. 
By shrinking the open intervals $I_j$ in the subcovering a bit we can assume 
that $K$ is covered by finitely many compact intervals $K_i$ each of which 
lies in some $I_j$. 
Since the local differentiable lifts have locally bounded 
derivatives, there exist constants $C_{K_i}$  for all $i$ such that
\[
\|\bar{c}'(t)\| \le C_{K_i} \qquad \mbox{for all}~ t \in K_i.
\] 
If we put $C_K := \max\{C_{K_i}\}$, then  
\[
\|\bar{c}'(t)\| \le C_{K} \qquad \mbox{for all}~ t \in K.
\] 
This completes the proof.
\endproof

\section{$C^1$-lifts} \label{secc1}

The following proposition is a slight modification of theorem 4.2 in 
\cite{rep-lift3}.

\begin{proposition}  \label{conder}
Assume that $\rho : G \to O(V)$ is a representation of a finite group 
$G$ with property $(\B_k)$.
Let
$c : \mathbb{R} \rightarrow V/G
= \sigma(V) \subseteq \mathbb{R}^n$ be a curve of class $C^{k+d}$. 
Then for any $t_0 \in \mathbb{R}$ there exists a local differentiable lift
$\bar{c}$ of $c$ near $t_0$ whose derivative is continuous at $t_0$.
\end{proposition}

\proof
Without loss of generality we may assume that $t_0 = 0$.
We show the existence of local differentiable lifts of $c$ whose derivatives 
are continuous at $0$ through any $v \in \sigma^{-1}(c(0))$.
By lemma \ref{fix}
we can assume $V^G =\{0\}$.

If $c(0) \ne 0$ corresponds to a regular orbit,
then unique orthogonal $C^{k+d}$-lifts defined near $0$ exist through
all $v \in \sigma^{-1}(c(0))$,
by lemma \ref{recall}(a).

If $c(0) = 0$,
then $c_1$ must vanish of at least second order at $0$,
since $c_1(t) \ge 0$ for all $t \in \mathbb{R}$.
That means $c_1(t) = t^2 c_{1,1}(t)$ near $0$ for a $C^{k+d-2}$-function
$c_{1,1}$.
By the multiplicity lemma \ref{mult}
we find that $c_i(t) = t^{d_i} c_{i,i}(t)$ near $0$ for $1 \le i \le n$,
where $c_{1,1},c_{2,2},\ldots,c_{n,n}$ are functions of class 
$C^{k+d-2},C^{k+d-d_2},\ldots,C^{k+d-d_n}$, respectively.
We consider the following $C^k$-curve in $\sigma(V)$ (since $\sigma(V)$ 
is closed in $\mathbb{R}^n$,
see \cite{procesischwarz}):
\begin{eqnarray*}
c_{(1)}(t) &:=& (c_{1,1}(t),
c_{2,2}(t),\ldots,c_{n,n}(t))\\ &=&
(t^{-2} c_1(t),t^{-d_2}c_2(t),\ldots,t^{-d_n} c_n(t)).
\end{eqnarray*}
By property $(\B_k)$ and proposition \ref{Beverywhere},
there exists a local differentiable lift $\bar{c}_{(1)}$ of $c_{(1)}$ 
with locally bounded derivative.
Thus, $\bar{c}(t) := t \cdot \bar{c}_{(1)}(t)$ is a 
local differentiable lift of $c$ near $0$ with derivative 
$\bar{c}'(t) = \bar{c}_{(1)}(t) + t \bar{c}_{(1)}'(t)$ 
which is continuous at $t =0$ with $\bar{c}'(0) = \bar{c}_{(1)}(0)$.
Note that $\sigma^{-1}(0) = \{0\}$,
therefore we are done in this case.

If $c(0) \ne 0$ corresponds to a singular orbit,
let $v$ be in $\sigma^{-1}(c(0))$ and consider the slice representation
$G_v \rightarrow O(N_v)$.
By \ref{slice},
the lifting problem reduces to the same problem for curves in $N_v/G_v$
now passing through $0$. 
By proposition \ref{Bk&slice} we may refer to the former case.
\endproof

\begin{theorem} \label{c1}
Assume that $\rho : G \to O(V)$ is a representation of a finite group 
$G$ with property $(\B_k)$.
Let
$c : \mathbb{R} \rightarrow V/G
= \sigma(V) \subseteq \mathbb{R}^n$ be a curve of class $C^{k+d}$. 
Then any differentiable lift $\bar{c}$ of $c$ is actually of class $C^1$.
\end{theorem}

\proof
Let $\bar{c}$ be a differentiable lift of $c$. 
Let $t_0 \in \mathbb{R}$ be arbitrary. 
We show that $\bar{c}'$ is continuous at $t_0$. 
Let $\tilde{c}$ denote the local differentiable lift of $c$ near $t_0$ 
with continuous derivative at $t_0$, provided by proposition \ref{conder}.
Consider a sequence $(t_m)_m \subseteq \mathbb{R}$ with $t_m \to t_0$. 
For every $m$ there is a $g_m \in G$ such that 
$\bar{c}(t_m) = g_m.\tilde{c}(t_m)$. 
Since $G$ is finite, we may choose a subsequence of $(t_m)_m$ 
again denoted by $(t_m)_m$ such that $\bar{c}(t_m) = g.\tilde{c}(t_m)$ 
for some fixed $g \in G$ and all $m$. 
By lemma \ref{deroflifts}, there exist $h_m \in G_{\bar{c}(t_m)}$ 
with $\bar{c}'(t_m) = h_m g.\tilde{c}'(t_m)$ for all $m$. 
Passing again to a subsequence we find a fixed $h \in G_{\bar{c}(t_m)}$ 
such that $\bar{c}(t_m) = h.\bar{c}(t_m) = h g.\tilde{c}(t_m)$ 
and $\bar{c}'(t_m) = h g.\tilde{c}'(t_m)$ for all $m$. 
Then
\[\bar{c}(t_0) = \lim_{t_m \to t_0} \bar{c}(t_m) 
  = \lim_{t_m \to t_0} h g.\tilde{c}(t_m) 
  = h g.\lim_{t_m \to t_0} \tilde{c}(t_m) = h g.\tilde{c}(t_0)\]
and
\[\bar{c}'(t_0) 
  = \lim_{t_m \to t_0} \frac{\bar{c}(t_m)-\bar{c}(t_0)}{t_m-t_0} 
  = \lim_{t_m \to t_0} \frac{h g.\tilde{c}(t_m)-h g.\tilde{c}(t_0)}{t_m-t_0} 
  = h g.\tilde{c}'(t_0)\] 
and hence
\[\lim_{t_m \to t_0} \bar{c}'(t_m) 
  = \lim_{t_m \to t_0} h g.\tilde{c}'(t_m) 
  = h g.\tilde{c}'(t_0) = \bar{c}'(t_0).\]
This completes the proof.
\endproof

The forgoing theorem \ref{c1} is false, if $G$ is not finite: 
\begin{example*}
Again consider the standard action of $SO(2)$ on $\mathbb{R}^2$ 
with orbit map $\sigma(x_1,x_2) = x_1^2+x_2^2$. 
Let us consider the curve $c(t) = t^4$ and its differentiable lift 
$\bar{c}(t) = \left(t^2 \cos \frac{1}{t},t^2 \sin \frac{1}{t}\right)$. 
But the derivative 
$\bar{c}'(t) = \left(2 t \cos \frac{1}{t} 
  + \sin \frac{1}{t},2 t \sin \frac{1}{t} - \cos \frac{1}{t}\right)$ 
is not continuous at $t = 0$. 
\end{example*}

\begin{remark*}
The failure of theorem \ref{c1} in this special example really 
is due to the fact that $SO(2)$ is infinite, 
since there is the following result due to Bony \cite{bony}: 
Any non-negative function $f : \mathbb{R} \to \mathbb{R}$ of class $C^{2m}$ 
can be represented as sum of squares $f = g^2 + h^2$ 
of $C^m$-functions $g$ and $h$. 
This result implies that $SO(2) : \mathbb{R}^2$ has property $(\B_2)$, 
and hence any standard representation of $SO(n)$ on $\mathbb{R}^n$ 
($n \ge 2$) has property $(\B_2)$ as well.
But see \ref{expolar}.

Note that the lifting problem for the standard representation of 
$SO(n)$ on $\mathbb R^n$ is just the problem of representing  
non-negative functions as sums of squares. 
In this regard consult \cite{bony0}, \cite{bony}, \cite{BBCP}, 
\cite{feffermanphong}, and \cite{hilbert}. 
\end{remark*}

\section{Twice differentiable lifts} \label{sec2der}

The following proposition modifies proposition \ref{conder}.

\begin{proposition}  \label{2der}
Assume that $\rho : G \to O(V)$ is a representation of a finite group 
$G$ with property $(\B_k)$.
Let
$c = (c_1,\ldots,c_n) : \mathbb{R} \rightarrow V/G
= \sigma(V) \subseteq \mathbb{R}^n$ be a curve of class $C^{k+2d}$.
Then for any $t_0 \in \mathbb{R}$ there exists a local $C^1$-lift
$\bar{c}$ of $c$ near $t_0$ which is twice differentiable at $t_0$.
\end{proposition}

\proof
Without loss of generality we may assume that $t_0 = 0$.
We show the existence of local $C^1$-lifts of $c$ which are 
twice differentiable at $0$ through any $v \in \sigma^{-1}(c(0))$.
By lemma \ref{fix}
we can assume $V^G =\{0\}$.

If $c(0) \ne 0$ corresponds to a regular orbit,
then unique orthogonal $C^{k+2d}$-lifts defined near $0$ exist through
all $v \in \sigma^{-1}(c(0))$,
by lemma \ref{recall}(a).

If $c(0) = 0$,
then as in the proof of proposition \ref{conder} we find that the curve
\begin{eqnarray*}
c_{(1)}(t) &:=& (c_{1,1}(t),
c_{2,2}(t),\ldots,c_{n,n}(t))\\ &=&
(t^{-2} c_1(t),t^{-d_2}c_2(t),\ldots,t^{-d_n} c_n(t))
\end{eqnarray*}
lies in $\sigma(V)$ and is of class $C^{k+d}$.
By property $(\B_k)$ and theorem \ref{c1},
there exists a local $C^1$-lift $\bar{c}_{(1)}$ of $c_{(1)}$.
Thus, $\bar{c}(t) := t \cdot \bar{c}_{(1)}(t)$ is a local $C^1$-lift of $c$
near $0$ with derivative 
$\bar{c}'(t) = \bar{c}_{(1)}(t) + t \bar{c}_{(1)}'(t)$ 
which is differentiable at $t=0$:
\[\lim_{t \to 0} \frac{\bar{c}'(t)-\bar{c}'(0)}{t} 
  = \lim_{t \to 0} \frac{\bar{c}_{(1)}(t) 
  - \bar{c}_{(1)}(0) + t \bar{c}_{(1)}'(t)}{t} = 2 \bar{c}_{(1)}'(0).\]
Note that $\sigma^{-1}(0) = \{0\}$,
therefore we are done in this case.

If $c(0) \ne 0$ corresponds to a singular orbit,
let $v$ be in $\sigma^{-1}(c(0))$ and consider the isotropy representation
$G_v \rightarrow O(N_v)$.
By \ref{slice},
the lifting problem reduces to the same problem for curves in $N_v/G_v$
now passing through $0$.
By proposition \ref{Bk&slice} we may refer to the former case.
\endproof

\begin{theorem} \label{twice}
Assume that $\rho : G \to O(V)$ is a representation of a 
finite group $G$ with property $(\B_k)$.
Let
$c : \mathbb{R} \rightarrow V/G
= \sigma(V) \subseteq \mathbb{R}^n$ be a curve of class $C^{k+2d}$. 
Then there exists a global twice differentiable lift $\bar{c}$ of $c$.
\end{theorem}

\proof
The proof will be carried out by induction on the cardinality of $G$.

If $G = \{e\}$ is trivial,
then $\bar{c} := c$ is a global twice differentiable lift.

So let us assume that for any finite $G'$ with $|G'| < |G|$ and any
$c : \mathbb{R} \rightarrow V/G'$
of class $C^{k+2d'}$ 
there exists a global twice differentiable lift
$\bar{c} : \mathbb{R} \rightarrow V$ of $c$,
where $\rho' : G' \rightarrow O(V)$ is an orthogonal representation on an
arbitrary real finite dimensional Euclidean vector space $V$ with property 
$(\B_k)$, and $d' = d(\rho')$.

We shall prove that the same is true for $G$.
Let
$c = (c_1,\ldots,c_n) : \mathbb{R} \rightarrow V/G
  = \sigma(V) \subseteq \mathbb{R}^n$ be of class $C^{k+2d}$.
We may assume that $V^G = \{0\}$,
by lemma \ref{fix}.
We can write
$c^{-1}(\sigma(V) \backslash \{0\}) = \bigcup_{i} (a_i,b_i)$,
a disjoint, at most countable union,
where $a_i,
b_i \in \mathbb{R} \cup \{\pm \infty\}$ with
$a_i < b_i$
such that each $(a_i,b_i)$ is maximal with respect to not containing zeros of $c$.
In particular,
we have $c(a_i) = c(b_i) = 0$ for all $a_i,
b_i \in \mathbb{R}$ appearing in the above presentation.

{\it Claim: On each $(a_i,b_i)$ there exists a twice differentiable lift
$\bar{c} : (a_i,b_i) \rightarrow V \backslash \{0\}$ of the restriction
$c|_{(a_i,b_i)} : (a_i,b_i) \rightarrow \sigma(V) \backslash \{0\}$}.
The lack of nontrivial fixed points guarantees that for all
$v \in V \backslash \{0\}$ the isotropy groups $G_v$ satisfy $|G_v| < |G|$.
Therefore, by induction hypothesis, which is fulfilled by proposition 
\ref{Bk&slice} and lemma \ref{dD}, and by \ref{slice},
we find local twice differentiable lifts of $c|_{(a_i,b_i)}$ near any
$t \in (a_i,b_i)$ and through all $v \in \sigma^{-1}(c(t))$.
Suppose that
$\bar{c}_1 : (a_i,b_i) \supseteq (a,b) \rightarrow V \backslash \{0\}$
is a local twice differentiable lift of $c|_{(a_i,b_i)}$ 
with maximal domain $(a,b)$, where,
say, $b < b_i$.
Then there exists a local twice differentiable lift $\bar{c}_2$ of
$c|_{(a_i,b_i)}$ near $b$,
and there exists a $t_0 < b$ such that both $\bar{c}_1$ and
$\bar{c}_2$ are defined near $t_0$.
Let $(t_m)_m$ be a sequence with $t_m \to t_0$. 
For any $m$ there exists a $g_m \in G$ such that 
$\bar{c}_1(t_m) = g_m.\bar{c}_2(t_m)$.
Since $G$ is finite, we may choose a subsequence again denoted by 
$(t_m)_m$ such that $\bar{c}_1(t_m) = g.\bar{c}_2(t_m)$ 
for a fixed $g \in G$ and for all $m$.
By lemma \ref{deroflifts}, there are $h_m \in G_{\bar{c}_1(t_m)}$ 
with $\bar{c}_1'(t_m) = h_m g.\bar{c}_2'(t_m)$ for all $m$. 
Passing again to a subsequence we find a fixed $h \in G_{\bar{c}_1(t_m)}$ 
such that $\bar{c}_1(t_m) = g.\bar{c}_2(t_m) = hg.\bar{c}_2(t_m)$ 
and $\bar{c}_1'(t_m) = h g.\bar{c}_2'(t_m)$ for all $m$. 
Consequently,
\[\bar{c}_1(t_0) = \lim_{t_m \to t_0} \bar{c}_1(t_m) 
= \lim_{t_m \to t_0} hg.\bar{c}_2(t_m) = hg.\bar{c}_2(t_0)\]
and
\[\bar{c}_1'(t_0) = \lim_{t_m \to t_0} \bar{c}_1'(t_m) 
  = \lim_{t_m \to t_0} hg.\bar{c}_2'(t_m) = hg.\bar{c}_2'(t_0),\]
and hence
\[\bar{c}_1''(t_0) = \lim_{t_m \to t_0} 
  \frac{\bar{c}_1'(t_m)-\bar{c}_1'(t_0)}{t_m-t_0} 
  = \lim_{t_m \to t_0} \frac{hg.\bar{c}_2'(t_m)-hg.\bar{c}_2'(t_0)}{t_m-t_0} 
  = hg.\bar{c}_2''(t_0).\]
So $\bar{c}(t) := \bar{c}_1(t)$ for $t \le t_0$ and 
$\bar{c}(t) := hg.\bar{c}_2(t)$ for $t \ge t_0$ defines a twice 
differentiable lift of $c|_{(a_i,b_i)}$ on a larger interval than $(a,b)$. 
This proves the claim.

Now let $\bar{c} : (a_i,b_i) \rightarrow V \backslash \{0\}$
be the twice differentiable lift of $c|_{(a_i,b_i)}$ constructed above.
For $a_i \ne - \infty$, we put $\bar{c}(a_i) := 0$ and 
$\bar{c}'(a_i) := \lim_{t \searrow a_i} \frac{\bar{c}(t)}{t-a_i}$ 
which exists as shown in the proof of theorem 4.4 in \cite{rep-lift3}. 
Then $\bar{c}$ is one-sided continuous at $a_i$, 
since $\langle \bar{c}(t) | \bar{c}(t) \rangle = \sigma_1(\bar{c}(t)) = c_1(t)$. 
Let $\tilde{c}$ be a local $C^1$-lift of $c$ defined near $a_i$ 
which is twice differentiable at $a_i$, provided by proposition \ref{2der}.
Then we find
\[\lim_{t \searrow a_i} \bar{c}(t) = \bar{c}(a_i) = 0 = \tilde{c}(a_i).\] 
Let $(t_m)_m \subseteq (a_i,b_i)$ be a sequence with $t_m \searrow a_i$. 
For any $m$ there exists a $g_m \in G$ such that 
$\bar{c}(t_m) = g_m.\tilde{c}(t_m)$.
We may choose a subsequence (again denoted by $(t_m)_m$) such that 
$\bar{c}(t_m) = g.\tilde{c}(t_m)$ for a fixed $g \in G$ and for all $m$.
By lemma \ref{deroflifts}, there are $h_m \in G_{\bar{c}(t_m)}$ with 
$\bar{c}'(t_m) = h_m g.\tilde{c}'(t_m)$ for all $m$. 
Passing again to a subsequence we find a fixed $h \in G_{\bar{c}(t_m)}$ 
such that $\bar{c}(t_m) = g.\tilde{c}(t_m) = hg.\tilde{c}(t_m)$ and 
$\bar{c}'(t_m) = h g.\tilde{c}'(t_m)$ for all $m$. 
Therefore we have
\begin{equation} \label{eqderbartilde}
\bar{c}'(a_i) = \lim_{t_m \searrow a_i} \frac{\bar{c}(t_m)}{t_m-a_i} 
  = \lim_{t_m \searrow a_i} \frac{hg.\tilde{c}(t_m)}{t_m-a_i} 
  = hg.\tilde{c}'(a_i).
\end{equation}
Moreover,
\[\lim_{t_m \searrow a_i} \bar{c}'(t_m) 
  = \lim_{t_m \searrow a_i} hg.\tilde{c}'(t_m) 
  = hg.\tilde{c}'(a_i) = \bar{c}'(a_i),\]
since $\tilde{c}$ is $C^1$. 
It follows that the set of all accumulation points of 
$(\bar{c}'(t))_{t \searrow a_i}$ lies in the orbit $G.\tilde{c}'(a_i)$. 
Since $G$ is finite, lemma \ref{recall}(b) implies that $\bar{c}'(t)$ 
converges for $t \searrow a_i$, with limit $\bar{c}'(a_i)$, 
because it does so along the sequence $(t_m)_m$. 
Otherwise put, the lift $\bar{c}$ is continuously differentiable 
also at the boundary point $a_i$ of its domain.

For the sequence $(t_m)_m$ from above we can argue further
\[\frac{\bar{c}'(t_m)-\bar{c}'(a_i)}{t_m-a_i} 
  = \frac{hg.\tilde{c}'(t_m)-hg.\tilde{c}'(a_i)}{t_m-a_i} 
  \to hg.\tilde{c}''(a_i) \qquad \mbox{as}~ t_m \searrow a_i, \]
since the lift $\tilde{c}$ is twice differentiable at $a_i$. 
Hence the set of all accumulation points of 
$\left(\frac{\bar{c}'(t)-\bar{c}'(a_i)}{t-a_i} \right)_{t \searrow a_i}$ 
is a subset of $G_{\bar{c}'(a_i)} hg.\tilde{c}''(a_i)$: 
Any accumulation point of 
$\left(\frac{\bar{c}'(t)-\bar{c}'(a_i)}{t-a_i} \right)_{t \searrow a_i}$ 
corresponds to a sequence 
$(t_m)_m \in (a_i,b_i)$ with $t_m \searrow a_i$ such that 
$\frac{\bar{c}'(t_m)-\bar{c}'(a_i)}{t_m-a_i} 
  \to \hat{h} \hat{g}.\tilde{c}''(a_i)$, 
where $\hat{h}$ and $\hat{g}$ are found by repeating the procedure above. 
From the equation 
$\hat{h} \hat{g}.\tilde{c}'(a_i) = \bar{c}'(a_i) = hg.\tilde{c}'(a_i)$, 
which follows from \eqref{eqderbartilde}, we can read off 
$(hg)^{-1} \hat{h} \hat{g} \in G_{\tilde{c}'(a_i)} 
  = (hg)^{-1} G_{\bar{c}'(a_i)} hg$, 
and hence $\hat{h} \hat{g} \in G_{\bar{c}'(a_i)} hg$.

By lemma \ref{recall}(b) we have that 
$\frac{\bar{c}'(t)-\bar{c}'(a_i)}{t-a_i}$ converges for 
$t \searrow a_i$, with limit $hg.\tilde{c}''(a_i)$, 
since it does so along the sequence $(t_m)_m$. 
That means that the one-sided second derivative of $\bar{c}$ exists at $a_i$. 
The same reasoning is true for $b_i \ne + \infty$. 
So we have extended our lift $\bar{c}$ twice differentiably 
to the closure of $(a_i,b_i)$.

Let us now construct a global twice differentiable lift of $c$ defined on the
whole of $\mathbb{R}$.
For isolated points $t_0 \in c^{-1}(0)$ the two twice differentiable lifts on
the neighboring intervals can be made to match twice differentiably, 
by applying a fixed transformation from $G$ to one of them: 
Let $\bar{c}_1$ and $\bar{c}_2$ denote the lifts left and right of $t_0$. 
Then $\bar{c}_1(t_0) = \bar{c}_2(t_0) = 0$ and, by lemma \ref{deroflifts}, 
we find some $g \in G$ such that $\bar{c}_1'(t_0) = g.\bar{c}_2'(t_0)$. 
Let $\tilde{c}$ be the local $C^1$-lift near $t_0$ which is 
twice differentiable at $t_0$, provided by proposition \ref{2der}. 
By the same argumentation as in the previous paragraph we find 
$h_1,h_2 \in G$ such that 
\[h_1.\tilde{c}'(t_0) = \bar{c}_1'(t_0) = g.\bar{c}_2'(t_0) 
  = h_2.\tilde{c}'(t_0),\]
and for the one-sided second derivatives we have
\[\lim_{t \nearrow t_0} \frac{\bar{c}_1'(t) - \bar{c}_1'(t_0)}{t-t_0} 
  = h_1.\tilde{c}''(t_0) \quad 
  \mbox{and} \quad \lim_{t \searrow t_0} 
  \frac{g.\bar{c}_2'(t) - g.\bar{c}_2'(t_0)}{t-t_0} = h_2.\tilde{c}''(t_0). 
\]
It follows that there is a $h : = h_1 h_2^{-1} \in G_{\bar{c}_1'(t_0)}$ 
with $\bar{c}_1''(t_0) = hg.\bar{c}_2''(t_0)$, which shows the assertion.
  
Let $E$ be the set of accumulation points of $c^{-1}(0)$.
For connected components of $\mathbb{R} \backslash E$ 
we can proceed inductively
to obtain twice differentiable lifts on them.

Let $\hat{c} : \mathbb{R} \to V$ be a global $C^1$-lift 
of $c$ which exists by theorem \ref{recall} and theorem \ref{c1}. 
We define the following set 
\[F := \{t \in \mathbb{R} : \hat{c}(t) = \hat{c}'(t) = 0\}.\] 

Note that every
lift $\bar c$ of $c$ has to vanish on $E$
and is continuous there since
$\langle \bar c(t)|\bar c(t) \rangle=\sigma_1(\bar c(t))=c_1(t)$.
We also claim that any lift
$\bar c$ of $c$ is differentiable at any point
$t'\in E$ with derivative 0.
Namely, the difference quotient $t\mapsto \frac{\bar c(t)}{t-t'}$ is a lift
of the curve 
\[
c_{(1,t')}(t) := ((t-t')^{-d_1} c_1(t),
  \ldots,(t-t')^{-d_n} c_n(t))
\]
in $\sigma(V)$ which vanishes at $t'$ by the following argument:
Consider the local lift $\tilde{c}$ of $c$ near $t'$,
provided by proposition \ref{2der}.
Let $(t_m)_{m \in \mathbb{N}} \subseteq c^{-1}(0)$ be a sequence with
$t'\ne t_m \rightarrow t'$,
consisting exclusively of zeros of $c$.
Such a sequence always exists since $t'\in E$.
Then we have
\[\tilde{c}'(t')
= \lim_{t \rightarrow t'} \frac{\tilde{c}(t) - \tilde{c}(t')}{t-t'}
= \lim_{m \rightarrow \infty} \frac{\tilde{c}(t_m)}{t_m -t'} = 0. \]
Thus $c_{(1,t')}(t')=\lim_{t\to t'}\sigma(\frac{\tilde c(t)}{t-t'})
=\sigma(\tilde c'(t'))=0$.

In particular this shows that $E \subseteq F$. 
If we denote by $F'$ the accumulation points of $F$, then 
$E \subseteq F = (F \backslash F') \cup F' \subseteq c^{-1}(0)$.

Consider first some $t' \in F \backslash F'$, i.e., 
$t'$ is an isolated point of $F$. 
Then again we have a local twice differentiable lift for $t \ne t'$ 
(left and right of $t'$), 
since near $t'$ there are only points of $\mathbb{R} \backslash E$. 
Moreover, proposition \ref{2der} yields again a local $C^1$-lift 
near $t'$ which is twice differentiable at $t'$. 
As above we are able to find a twice differentiable lift on the set 
$(\mathbb{R} \backslash E) \cup (F \backslash F')$.

Finally let $t' \in F'$, i.e., $t'$ is an accumulation point of $F$. 
By proposition \ref{2der}, we have again a local $C^1$-lift $\tilde{c}$ 
near $t'$ which is twice differentiable at $t'$. 
Lemma \ref{deroflifts} implies that locally near $t'$ the set $F$ 
is given by $F = \{\tilde{c}(t) = \tilde{c}'(t) = 0\}$. 
So we have $\tilde{c}(t') = \tilde{c}'(t') = \tilde{c}''(t') = 0$, 
as $t'$ is an accumulation point of $F$. 
We extend our twice differentiable lift $\bar{c}$ on 
$(\mathbb{R} \backslash E) \cup (F \backslash F')$ by $0$ on $F'$ 
to the whole of 
$(\mathbb{R} \backslash E) \cup (F \backslash F') \cup F' 
  = (\mathbb{R} \backslash E) \cup F = \mathbb{R}$. 
It remains to check that then $\bar{c}$ is twice differentiable at $t' \in F'$. 
Since $F' \subseteq E$, we obtain that $\bar{c}$ vanishes at $t'$ and 
is continuous and differentiable there with derivative $0$. 
Consider a sequence $(t_m)_m$ with $t' \ne t_m \to t'$. 
Passing to subsequences, we find, for all $m$,  
$\bar{c}(t_m) = g.\tilde{c}(t_m)$ and $\bar{c}'(t_m) = hg.\tilde{c}'(t_m)$ 
for some $g \in G$ and some $h \in G_{\bar{c}(t_m)}$, 
by lemma \ref{deroflifts}. 
Then,
\[ \frac{\bar{c}'(t_m) - \bar{c}'(t')}{t_m - t'} 
  =  \frac{\bar{c}'(t_m)}{t_m - t'} 
  =  \frac{hg.\tilde{c}'(t_m)}{t_m - t'}  \to hg.\tilde{c}''(t') 
  = 0 \qquad \mbox{as} ~ t_m \to t'.\]
It follows that the second derivative of $\bar{c}$ at $t'$ 
exists and equals $0$. 
This completes the proof. 
\endproof

\section{Polar representations} \label{secpolar}

The main results of sections \ref{secc1} and \ref{sec2der}, 
obtained there for finite groups $G$, can be generalized 
to polar representations $G \to O(V)$.  

An orthogonal representation $\rho : G \rightarrow O(V)$ of a Lie group $G$ 
on a finite dimensional Euclidean vector space $V$ is called \emph{polar}, 
if there exists a linear subspace $\Sigma \subseteq V$, 
called a \emph{section} or a \emph{Cartan subspace}, 
which meets each orbit orthogonally. 
See \cite{Dadok}, \cite{DadokKac}, and \cite{PalaisTerng0}. 
The trace of the $G$-action is the action of the \emph{generalized Weyl group} 
$W(\Sigma) = N_G(\Sigma)/Z_G(\Sigma)$ on $\Sigma$, where 
$N_G(\Sigma) := \{g \in G : \rho(g)(\Sigma) = \Sigma\}$ and 
$Z_G(\Sigma) := \{g \in G : \rho(g)(s) = s ~\mbox{for all}~ s \in \Sigma\}$. 
The generalized Weyl group is a finite group, and is a reflection group if 
$G$ is connected. If $\Sigma'$ is a different section, 
then there is an isomorphism $W(\Sigma) \to W(\Sigma')$ 
induced by an inner automorphism of $G$.

We shall need the following generalization of Chevalley's restriction theorem, 
which is due to Dadok and Kac \cite{DadokKac} 
and independently to Terng \cite{Terng}.

\begin{theorem} \label{restriction}
Let $\rho : G \rightarrow O(V)$ be a polar representation 
of a compact Lie group, with section $\Sigma$ and 
generalized Weyl group $W(\Sigma)$. 
Then the algebra $\mathbb{R}[V]^G$ of $G$-invariant polynomials on $V$ 
is isomorphic to the algebra $\mathbb{R}[\Sigma]^{W(\Sigma)}$ 
of $W(\Sigma)$-invariant polynomials on the section $\Sigma$, 
via restriction $f \mapsto f|_{\Sigma}$.
\end{theorem}

As a consequence of this theorem we obtain that the orbit spaces 
$V/G = \sigma(V)$ and $\Sigma/W(\Sigma) = \sigma |_{\Sigma}(\Sigma)$ 
are isomorphic.

\begin{theorem} \label{polar}
Let $\rho : G \rightarrow O(V)$ be a polar representation of a 
compact Lie group on a finite dimensional Euclidean vector space $V$ 
with orbit map $\sigma : V \rightarrow \mathbb{R}^n$. 
Assume that $W(\Sigma) \to O(\Sigma)$ has property $(\B_k)$ for some 
section $\Sigma$. 
Let $c : \mathbb{R} \rightarrow \sigma(V) \subseteq \mathbb{R}^n$ 
be a curve in the orbit space. 
Then we have:
\begin{enumerate}
\item[(1)] If $c$ is of class $C^{k+d}$, 
  then there exists a global orthogonal 
  $C^1$-lift $\bar{c} : \mathbb{R} \to V$.
\item[(2)] If $c$ is of class $C^{k+2d}$, 
  then there exists a global orthogonal twice differentiable lift 
  $\bar{c} : \mathbb{R} \to V$. 
\end{enumerate}
\end{theorem}

\proof
By theorem \ref{restriction}, $\sigma |_{\Sigma} : \Sigma \to \mathbb{R}^n$ 
is the orbit map for the representation $W(\Sigma) \to O(\Sigma)$, 
and hence the orbit spaces 
$V/G = \sigma(V)$ and $\Sigma/W(\Sigma) = \sigma |_{\Sigma}(\Sigma)$ 
are isomorphic.

If $c : \mathbb{R} \to \sigma(V) \cong \sigma |_{\Sigma}(\Sigma)$ is 
$C^{k+d}$, then by theorem \ref{recall} and theorem 
\ref{c1} (since $W(\Sigma)$ is finite) there exists a global $C^1$-lift 
$\bar{c} : \mathbb{R} \to \Sigma$, which as a curve in $V$ is orthogonal 
to each $G$-orbit it meets, by the properties of $\Sigma$. 
This shows $(1)$.

If $c : \mathbb{R} \to \sigma(V) \cong \sigma |_{\Sigma}(\Sigma)$ 
is $C^{k+2d}$, then statement $(2)$ follows analogously from theorem \ref{2der}.
\endproof

\begin{example} \label{expolar}
The standard representation of $SO(n)$ on $\mathbb{R}^n$ is polar. 
Any $1$-dimensional linear subspace $\Sigma$ of $\mathbb{R}^n$ is a section. 
The associated generalized Weyl group is $W(\Sigma) = \{\pm \id\}$. 
So the representation $W(\Sigma) \to O(\Sigma)$ has property $(\B_2)$, 
since it reduces to $S_2 : \mathbb R^2$, 
the problem of finding regular roots of $x^2 - f(t) = 0$ ($f \ge 0$ and $C^2$), 
which has property $(\B_2)$, by theorem \ref{bron-wa}. 
Hence theorem \ref{polar} is applicable.
\end{example}

\section{Stability of property $(\B)$} \label{secstab}

\begin{proposition} \label{subrep}
Let $\rho : G \to O(V)$ be an orthogonal representation of a 
compact Lie group on a real finite dimensional Euclidean vector space 
$V$ having property $(\B_k)$. 
For any $G$-invariant linear subspace $W \subseteq V$ 
the subrepresentation $\rho' : G \to O(W)$ has property $(\B_k)$ as well. 
%Note that $k \ge d \ge d'$, where $d$ and $d'$ are the integers 
%associated to the representations $\rho$ and $\rho'$ as in \ref{setting}.
\end{proposition}

\proof
Let us consider the following restriction map
\begin{eqnarray*}
\mathbb{R}[V]^G &\longrightarrow& \mathbb{R}[W]^G\\ p &\longmapsto& p|_W.
\end{eqnarray*}
This map is a surjective algebra homomorphism.
%, and it is surjective: 
%For $q \in \mathbb{R}[W]^G$ and the $G$-invariant decomposition 
%$V = W \oplus W^{\bot}$ we define $p(w_1 + w_2) := q(w_1)$. 
%Then $p \in \mathbb{R}[V]^G$ and $p|_W = q$.
%
Hence, if $\sigma_1,\ldots,\sigma_n$ are generators of $\mathbb{R}[V]^G$, 
then their restrictions $\sigma_1|_W,\ldots,\sigma_n|_W$ 
generate $\mathbb{R}[W]^G$:
\[\mathbb{R}[W]^G = \mathbb{R}\left[\sigma_1|_W,\ldots,\sigma_n|_W\right].\]
Of course, $\sigma_1|_W,\ldots,\sigma_n|_W$ may not be a 
minimal system of generators of $\mathbb{R}[W]^G$. 
Nevertheless, the map 
$\sigma|_W = (\sigma_1|_W,\ldots,\sigma_n|_W) : 
  W \to \sigma |_W(W) \subseteq \mathbb{R}^n$ 
still induces a homeomorphism between the orbit space $W/G$ 
and the image $\sigma |_W(W)$, by \ref{setting}.
The following commuting diagram
\[
\xymatrix{
V \ar[r]^{\sigma} & \sigma(V)\\
W \ar@{^{(}->}[u] \ar[r]^{\sigma |_W} & \sigma|_W(W) \ar@{^{(}->}[u] 
}
\]
then shows that the orbit space $W/G = \sigma|_W(W)$ is a 
(in general lower dimensional) subset of the orbit space $V/G = \sigma(V)$.

Let $c : \mathbb{R} \to \sigma|_W(W) \cap U$ be a $C^{k}$-curve in the 
orbit space $\sigma|_W(W)$, where $U = U(\rho)$ is the open neighborhood of 
$0$ in $\sigma(V)$ 
from the definition of property $(\B_k)$ (see \ref{propertyB}). 
We may view $c$ as a curve in the orbit space $\sigma(V)$, 
and since the representation $\rho$ has property $(\B_k)$, 
we can lift $c$ to a local differentiable curve $\bar{c}$ in $V$ 
with locally bounded derivative. 
But then $\bar{c}$ has obviously to lie in the $G$-invariant subspace $W$.
%For all other slice representations of $\rho'$ we can argue in the same way: 
%For $w \in W$ we have $T_w(G.w) \subseteq W$, since $W$ is $G$-invariant, 
%and therefore
%\[T_w(G.w)^{\bot} = T_w(G.w)^{\bot_W} \oplus W^{\bot},\]
%where $\bot_W$ indicates the orthogonal complement with respect 
%to the inner product $\langle \quad | \quad \rangle$ restricted to $W$ in $W$. 
This completes the proof.
\endproof

\begin{proposition} \label{directsum}
Suppose that $\rho_i : G_i \to O(V_i)$ $(1 \le i \le l)$ are 
orthogonal representations of compact Lie groups $G_i$ on real 
finite dimensional Euclidean vector spaces $V_i$ having property $(\B_{k_i})$. 
Then the orthogonal direct sum 
\[\rho_1 \oplus \cdots \oplus \rho_l : 
  G_1 \times \cdots \times G_l 
  \longrightarrow O(V_1 \oplus \cdots \oplus V_l)\]
of the representations $\rho_1,\ldots,\rho_l$ has property $(\B_k)$, 
where $k = \max\{k_1,\ldots,k_l\}$. 
\end{proposition}

\proof
It is sufficient to consider the case $l=2$, since the general case 
follows by induction.

If $\langle \quad|\quad \rangle_1$ and $\langle \quad|\quad \rangle_2$ 
denote the inner products on $V_1$ and $V_2$, then 
\[\langle v_1 + v_2 | w_1 + w_2 \rangle 
  := \langle v_1 | w_1 \rangle_1 + \langle v_2 | w_2 \rangle_2  \]
defines an inner product on $V = V_1 \oplus V_2$ which makes 
$V_1$ and $V_2$ into orthogonal subspaces of $V$. 
The action of $G = G_1 \times G_2$ on $V_1 \oplus V_2$ 
is obviously again orthogonal. 
Moreover, we find 
$\mathbb{R}[V]^G 
  = \mathbb{R}[V_1 \oplus V_2]^{G_1 \times G_2} \cong 
  \mathbb{R}[V_1]^{G_1} \otimes \mathbb{R}[V_2]^{G_2}$ 
and $V/G = (V_1 \oplus V_2)/(G_1 \times G_2) \cong V_1/G_1 \times V_2/G_2$. 
%Clearly, for the isotropy groups one has 
%$(G_1 \times G_2)_{v_1 +v_2} = (G_1)_{v_1} \times (G_2)_{v_2}$. 
%Consider 
%\begin{eqnarray*}
%T_{v_1+v_2}((G_1 \times G_2).(v_1+v_2)) &\cong& 
%  T_{(v_1,v_2)}(G_1.v_1 \times G_2.v_2)\\ 
%  &\cong& T_{v_1}(G_1.v_1) \oplus T_{v_2}(G_2.v_2).
%\end{eqnarray*}
%Hence
%\begin{eqnarray*}
%N_{v_1+v_2} &=& T_{v_1+v_2}((G_1 \times G_2).(v_1+v_2))^{\bot} 
%  = (T_{v_1}(G_1.v_1) \oplus T_{v_2}(G_2.v_2))^{\bot}\\ 
%  &=& T_{v_1}(G_1.v_1)^{\bot} \cap T_{v_2}(G_2.v_2)^{\bot} \\ 
%  &=& (T_{v_1}(G_1.v_1)^{\bot_1} \oplus V_2) 
%    \cap(V_1 \oplus T_{v_2}(G_2.v_2)^{\bot_2})\\ 
%  &=& (N_1)_{v_1} \oplus (N_2)_{v_2},
%\end{eqnarray*}
%where $(N_i)_{v_i} = T_{v_i}(G_i.v_i)^{\bot_i}$ is the orthogonal complement 
%of $T_{v_i}(G_i.v_i)$ in $V_i$ 
%with respect to $\langle \quad | \quad \rangle_i$. 
%It follows 
%$\mathbb{R}[N_{v_1+v_2}]^{G_{v_1+v_2}} \cong 
%  \mathbb{R}[(N_1)_{v_1}]^{(G_1)_{v_1}} \otimes 
%  \mathbb{R}[(N_2)_{v_2}]^{(G_2)_{v_2}}$ 
%and 
%$N_{v_1+v_2}/G_{v_1+v_2} 
%  \cong (N_1)_{v_1}/(G_1)_{v_1} \times (N_2)_{v_2}/(G_2)_{v_2}$.
%
%Let $d_1$, $d_2$, and $d$ be the integers associated to 
%representations in \ref{setting} of $\rho_1$, $\rho_2$, 
%and $\rho = \rho_1 \oplus \rho_2$. 
%Then $d = \max\{d_1,d_2\}$, and therefore $k \ge d$. 

Now any $C^k$-curve $c$ in $U_1 \times U_2 \subseteq V/G$ has the form 
$c = (c_1,c_2)$ for $C^k$-curves $c_i$ in $U_i \subseteq V_i/G_i$, where 
$U_i = U(\rho_i)$ from \ref{propertyB}, 
which allow local differentiable lifts $\bar{c}_i$ with 
locally bounded derivative to $V_i$, by assumption. 
This shows that $\rho = \rho_1 \oplus \rho_2$ has property $(\B_k)$.
\endproof

%\begin{corollary} 
%A representation $\rho$ of a compact Lie group has
%property $(\B_k)$ if and only if each irreducible subrepresentation of $\rho$
%has property $(\B_k)$.
%\end{corollary}
%
%\proof This follows from \ref{subrep} for one direction, and from
%\ref{directsum} and \ref{subrep} for the converse.
%\endproof

\section{Finite groups $G$ have property $(\B)$} \label{secallfiniteG}

\begin{theorem} \label{finiteG}
Let $\rho : G \to O(V)$ be a real finite-dimensional 
orthogonal representation of a finite group $G$, 
and let $\sigma_1,\ldots,\sigma_n$ be a minimal system of 
homogeneous generators of $\mathbb{R}[V]^G$. 
Write $V = V_1 \oplus \cdots \oplus V_l$ as orthogonal direct sum of 
irreducible subspaces $V_i$. 
Choose $v_i \in V_i \backslash \{0\}$ such that the cardinality of the  
corresponding isotropy group $G_{v_i}$ is maximal, 
and put $k = \max\{d(\rho), |G|/|G_{v_i}| : 1 \le i \le l\}$.
Then any curve 
$c = (c_1,\ldots,c_n) : \mathbb{R} \to V/G = \sigma(V) \subseteq \mathbb{R}^n$ 
of class $C^{k}$ in the orbit space admits a global differentiable lift 
$\bar{c}$ to $V$ with locally bounded derivative. 
\end{theorem}

\proof
We shall reduce to the representation $S_{k} : \mathbb{R}^{k}$, 
the polynomial case.

Since $d \le k$, we can apply theorem \ref{recall} 
which provides a differentiable lift 
$\bar{c} : \mathbb{R} \to V$ 
of $c$.

Let $i$ be fixed. 
For $g \in G$ we define a linear function 
\begin{align*}
F_{i,g} : V &\longrightarrow \mathbb{R}\\
x &\longmapsto \langle v_i | g.\pr_{V_i}(x) \rangle = 
\langle v_i | g.x \rangle,
\end{align*}
since $V = V_1 \oplus \cdots \oplus V_l$ is an orthogonal direct sum. 
Here $\pr_{V_i} : V \to V_i$ is the natural projection.
The cardinality of distinct functions $F_{i,g}$ equals $k_i :=|G|/|G_{v_i}|$. 

Let $G_{v_i} \backslash G$ denote the space of right cosets of 
$G_{v_i}$ in $G$, and introduce a numbering $G_{v_i} \backslash G = \{g_1,g_2,\ldots,g_{k_i}\}$. 
We construct the following polynomials on $V$:
\begin{equation*} 
a_{i,j}(x) = \sum_{1 \le m_1 < \cdots < m_j \le k_i} F_{i,g_{m_1}}(x) \cdots F_{i,g_{m_j}}(x) 
\qquad 1 \le j \le k_i.
\end{equation*}

These polynomials $a_{i,j}$ are $G$-invariant by construction, and therefore 
expressible in the homogeneous generators $\sigma_1,\ldots,\sigma_n$ of $\mathbb{R}[V]^G$, i.e., there exist polynomials $p_{i,j} \in \mathbb{R}[\mathbb{R}^n]$ such that 
\begin{equation} \label{avssigma}
a_{i,j} = p_{i,j}(\sigma_1,\ldots,\sigma_n) 
  \qquad 1 \le j \le k_i.
\end{equation}
 
The polynomials $a_{i,j}$, for $1 \le j \le k_i$, 
are elementary symmetric functions in the variables $F_{i,g}(x)$, 
where $g$ runs through $G_{v_i} \backslash G = \{g_1,g_2,\ldots,g_{k_i}\}$.
Finally, we associate the following monic polynomial of degree 
$k_i$ in one variable $y$:
\begin{equation} \label{asspoly}
P_i(x)(y) = y^{k_i} + \sum_{j=1}^{k_i} (-1)^j a_{i,j}(x) y^{k_i-j} 
= \prod_{j=1}^{k_i} (y - F_{i,g_j}(x)).
\end{equation}
By construction, the functions $x \mapsto F_{i,g}(x)$ 
($g \in G_{v_i} \backslash G$) parameterize the roots of $x \mapsto P_i(x)(y)$ which, 
consequently, are always real. 

Now consider the functions $t \mapsto a_{i,j}(\bar{c}(t))$ 
($1 \le j \le k_i$) which are of class $C^{k}$ 
by equation \eqref{avssigma}. 
As in \eqref{asspoly} we may associate a $C^{k}$-curve 
$t \mapsto P_i(t)(y)$ of monic polynomials of degree 
$k_i$ in one variable defined by
\[P_i(t)(y) 
  = y^{k_i} + \sum_{j=1}^{k_i} (-1)^j a_{i,j}(\bar{c}(t)) y^{k_i-j}.\] 
%By construction, the differentiable functions $t \mapsto F_{i,g}(\bar{c}(t))$ 
%($g \in G_{v_i} \backslash G$) parameterize the roots of $t \mapsto P_i(t)(y)$ which, 
%consequently, are always real. 
By theorem \ref{bron-wa}, applied to the curve of polynomials 
$t \mapsto P_i(t)(y)$, 
the differentiable functions $t \mapsto F_{i,g}(\bar{c}(t))$ 
($g \in G_{v_i} \backslash G$) which parameterize the roots of 
$t \mapsto P_i(t)(y)$ have locally bounded derivative. 

Since $V_i$ is irreducible, the linear span of the orbit $G.v_i$ spans $V_i$. 
If we repeat the above procedure for each $1 \le i \le l$, it follows that 
$\bar{c}$ is a differentiable lift of $c$ with locally bounded derivative. 
This completes the proof.
\endproof

\begin{corollary} \label{cor}
Any real finite-dimensional orthogonal representation $\rho : G \to O(V)$ 
of a finite group $G$ has property $(\B_k)$, where 
$k = \max\{d(\rho), |G|/|G_{v_i}| : 1 \le i \le l\}$, 
$v_i \in V_i \backslash \{0\}$ are chosen such that the cardinality of 
$G_{v_i}$ is maximal, and $V = V_1 \oplus \cdots \oplus V_l$ is the 
decomposition into irreducible subrepresentations.\qed
\end{corollary} 

\begin{corollary}
Any polar representation $\rho$ of a compact Lie group $G$ 
has property $(\B_k)$, where $k$ is determined analogously to corollary 
\ref{cor} but for the representation $W \to O(\Sigma)$, where 
$W$ is the generalized Weyl 
group of some section $\Sigma$. 
Moreover, the lifts can be chosen orthogonal. \qed
\end{corollary}

\begin{remark*}
The case $k = |G|$ can occur: For finite rotation groups in the plane we have 
$d = |G|$, and for any non-zero $v$ the isotropy group $G_v$ is trivial.
\end{remark*}

\begin{remark*}
There are irreducible orthogonal representations of finite groups $G$ where
the inequality $d \le |G|/|G_v|$ is violated for non-zero vectors $v$:

Consider the rotational symmetry group $T$ of the regular tetrahedron in 
$\mathbb{R}^3$. We find $d=6$ (e.g. \cite{CumminsPatera}). 
The isotropy group of each vertex $v$ of the tetrahedron has $3$ elements. 
So $|G|/|G_v| =12/3=4$.  

Furthermore, the same phenomenon appears for the rotational symmetry groups 
$W$ and $H$ of the cube and the regular icosahedron in $\mathbb{R}^3$, 
respectively. Consult section \ref{rotgroups}.
\end{remark*}

\section{Property $(\B_k)$ for finite reflection groups} \label{reflectiongroups}

Abusing notation we will denote finite reflection groups 
as well as their root systems (respectively their Coxeter graphs) 
with the same symbols.

Recall the characterization of finite reflection groups 
(\cite{GroveBenson}, \cite{humphreys}):

If $G$ is a finite subgroup of $O(V)$ that is generated by reflections, 
then $V$ may be written as the orthogonal direct sum of 
$G$-invariant subspaces $V_0 = V^G,V_1,\ldots,V_k$ 
with the following properties:
\begin{enumerate}
\item[(a)] If $G_i = \{g|_{V_i} : g \in G\}$, 
  then $G_i$ is a subgroup of $O(V_i)$, and $G$ is isomorphic with 
  $G_0 \times G_1 \times \cdots \times G_k$.
\item[(b)] $G_0$ consists only of the identity transformation on $V_0$.
\item[(c)] Each $G_i$ ($i \ge 1$) is one of the groups 
\begin{gather*}
A_n, n\ge 1; B_n, n \ge 2; D_n, n \ge 4; I_2^n, n \ge 5, n \ne 6; \\ 
G_2; H_3; H_4; F_4; E_6; E_7; E_8.
\end{gather*}
\end{enumerate}

We will apply theorem \ref{finiteG} to each of the irreducible finite 
reflection groups listed in $(c)$. 
Here the inequality $d \le |G|/|G_v|$ will be satisfied for all non-zero $v$ 
which can be checked directly in the table given at the end of this section.
For irreducible finite reflection groups we may choose $v$ to be some 
non-zero vector in a one dimensional face of a Weyl chamber, in order to 
achieve that $|G|/|G_v|$ is minimal.

Let $\rho : G \to O(V)$ be the standard representation of some irreducible 
finite reflection group $G$ listed in $(c)$. 
We consider an arbitrary slice representation 
$G_v \to O(N_v)$ ($v \in V)$ of $\rho$; note $N_v = V$. 
By \cite{humphreys}, theorem 1.12, there exists a $g \in G$ such that 
$g.v = w$ for a $w$ in the fundamental domain 
$F = \{x \in V : \langle x | r \rangle \ge 0 ~\mbox{for all}~ r \in \Pi\}$,
where $\Pi$ is a system of simple roots, 
and $G_w$ is generated by the simple reflections it contains. 
It follows that we can read off easily the information we need to determine 
a minimal $|G|/|G_v|$ from the Coxeter graph of $G$.
Therefore, we give a complete list of all Coxeter graphs:

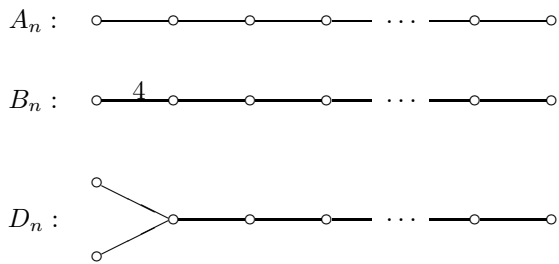
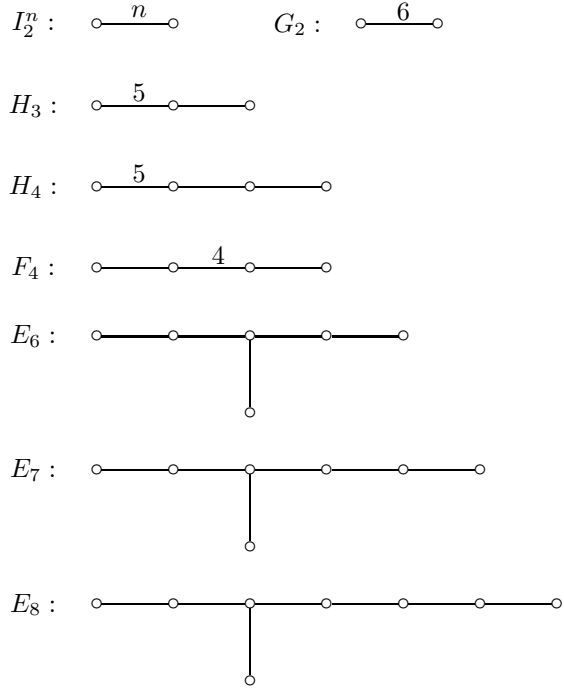
\[
\begin{picture}(190,15)
\put(-10,5){\makebox(0,0){$A_n:$}}
\put(14,5){\makebox(0,0){$\circ$}}
\put(16,5){\line(1,0){25}}
\put(43,5){\makebox(0,0){$\circ$}}
\put(45,5){\line(1,0){25}}
\put(72,5){\makebox(0,0){$\circ$}}
\put(74,5){\line(1,0){25}}
\put(101,5){\makebox(0,0){$\circ$}}
\put(103,5){\line(1,0){15}}
\put(130,5){\makebox(0,0){$\dots$}}
\put(155,5){\line(-1,0){15}}
\put(157,5){\makebox(0,0){$\circ$}}
\put(159,5){\line(1,0){25}}
\put(186,5){\makebox(0,0){$\circ$}}
\end{picture}
\]

\[
\begin{picture}(190,15)
\put(-10,5){\makebox(0,0){$B_n:$}}
\put(30,9){\makebox(0,0){4}}
\put(14,5){\makebox(0,0){$\circ$}}
\put(16,5){\line(1,0){25}}
\put(43,5){\makebox(0,0){$\circ$}}
\put(45,5){\line(1,0){25}}
\put(72,5){\makebox(0,0){$\circ$}}
\put(74,5){\line(1,0){25}}
\put(101,5){\makebox(0,0){$\circ$}}
\put(103,5){\line(1,0){15}}
\put(130,5){\makebox(0,0){$\dots$}}
\put(155,5){\line(-1,0){15}}
\put(157,5){\makebox(0,0){$\circ$}}
\put(159,5){\line(1,0){25}}
\put(186,5){\makebox(0,0){$\circ$}}
\end{picture}
\]     

\[
\begin{picture}(190,45)
\put(-10,20){\makebox(0,0){$D_n:$}}
\put(14,34){\makebox(0,0){$\circ$}}
\put(16,33){\line(2,-1){25}}
\put(14,6){\makebox(0,0){$\circ$}}
\put(16,7){\line(2,1){25}}
\put(43,20){\makebox(0,0){$\circ$}}
\put(45,20){\line(1,0){25}}
\put(72,20){\makebox(0,0){$\circ$}}
\put(74,20){\line(1,0){25}}
\put(101,20){\makebox(0,0){$\circ$}}
\put(103,20){\line(1,0){15}}
\put(130,20){\makebox(0,0){$\dots$}}
\put(155,20){\line(-1,0){15}}
\put(157,20){\makebox(0,0){$\circ$}}
\put(159,20){\line(1,0){25}}
\put(186,20){\makebox(0,0){$\circ$}}
\end{picture}
\]

\[
\begin{picture}(190,15)
\put(-10,5){\makebox(0,0){$I_2^n:$}}
\put(30,10){\makebox(0,0){$n$}}
\put(14,5){\makebox(0,0){$\circ$}}
\put(16,5){\line(1,0){25}}
\put(43,5){\makebox(0,0){$\circ$}}
\put(90,5){\makebox(0,0){$G_2:$}}
\put(130,10){\makebox(0,0){$6$}}
\put(114,5){\makebox(0,0){$\circ$}}
\put(116,5){\line(1,0){25}}
\put(143,5){\makebox(0,0){$\circ$}}
\end{picture}
\]

\[
\begin{picture}(190,15)
\put(-10,5){\makebox(0,0){$H_3:$}}
\put(30,10){\makebox(0,0){$5$}}
\put(14,5){\makebox(0,0){$\circ$}}
\put(16,5){\line(1,0){25}}
\put(43,5){\makebox(0,0){$\circ$}}
\put(45,5){\line(1,0){25}}
\put(72,5){\makebox(0,0){$\circ$}}
\end{picture}
\]

\[
\begin{picture}(190,15)
\put(-10,5){\makebox(0,0){$H_4:$}}
\put(30,10){\makebox(0,0){$5$}}
\put(14,5){\makebox(0,0){$\circ$}}
\put(16,5){\line(1,0){25}}
\put(43,5){\makebox(0,0){$\circ$}}
\put(45,5){\line(1,0){25}}
\put(72,5){\makebox(0,0){$\circ$}}
\put(74,5){\line(1,0){25}}
\put(101,5){\makebox(0,0){$\circ$}}
\end{picture}
\]

\[
\begin{picture}(190,15)
\put(-10,5){\makebox(0,0){$F_4:$}}
\put(60,10){\makebox(0,0){$4$}}
\put(14,5){\makebox(0,0){$\circ$}}
\put(16,5){\line(1,0){25}}
\put(43,5){\makebox(0,0){$\circ$}}
\put(45,5){\line(1,0){25}}
\put(72,5){\makebox(0,0){$\circ$}}
\put(74,5){\line(1,0){25}}
\put(101,5){\makebox(0,0){$\circ$}}
\end{picture}
\]

\[
\begin{picture}(190,35)
\put(-10,30){\makebox(0,0){$E_6:$}}
\put(14,30){\makebox(0,0){$\circ$}}
\put(16,30){\line(1,0){25}}
\put(43,30){\makebox(0,0){$\circ$}}
\put(45,30){\line(1,0){25}}
\put(72,30){\makebox(0,0){$\circ$}}
\put(74,30){\line(1,0){25}}
\put(101,30){\makebox(0,0){$\circ$}}
\put(103,30){\line(1,0){25}}
\put(130,30){\makebox(0,0){$\circ$}}
\put(72,28.5){\line(0,-1){25}}
\put(72,1){\makebox(0,0){$\circ$}}
\end{picture}
\]

\[
\begin{picture}(190,35)
\put(-10,30){\makebox(0,0){$E_7:$}}
\put(14,30){\makebox(0,0){$\circ$}}
\put(16,30){\line(1,0){25}}
\put(43,30){\makebox(0,0){$\circ$}}
\put(45,30){\line(1,0){25}}
\put(72,30){\makebox(0,0){$\circ$}}
\put(74,30){\line(1,0){25}}
\put(101,30){\makebox(0,0){$\circ$}}
\put(103,30){\line(1,0){25}}
\put(130,30){\makebox(0,0){$\circ$}}
\put(132,30){\line(1,0){25}}
\put(159,30){\makebox(0,0){$\circ$}}
\put(72,28.5){\line(0,-1){25}}
\put(72,1){\makebox(0,0){$\circ$}}
\end{picture}
\]

\[
\begin{picture}(190,35)
\put(-10,30){\makebox(0,0){$E_8:$}}
\put(14,30){\makebox(0,0){$\circ$}}
\put(16,30){\line(1,0){25}}
\put(43,30){\makebox(0,0){$\circ$}}
\put(45,30){\line(1,0){25}}
\put(72,30){\makebox(0,0){$\circ$}}
\put(74,30){\line(1,0){25}}
\put(101,30){\makebox(0,0){$\circ$}}
\put(103,30){\line(1,0){25}}
\put(130,30){\makebox(0,0){$\circ$}}
\put(132,30){\line(1,0){25}}
\put(159,30){\makebox(0,0){$\circ$}}
\put(161,30){\line (1,0){25}}
\put(188,30){\makebox(0,0){$\circ$}}
\put(72,28.5){\line(0,-1){25}}
\put(72,1){\makebox(0,0){$\circ$}}
\end{picture}
\]

Easy computations yield the results collected in the following 
table which gives a complete survey of the standard representations 
of all irreducible finite reflection groups. 
The integers $d$ and $|G|$ for the listed representations can be found 
e.g. in \cite{GroveBenson}; 
in \cite{mehta} also generators of the corresponding algebra of invariant 
polynomials are available. 
The integer $k$ is the 
minimum of the numbers $|G|/|G_v|$ where $v$ runs through all non-zero 
vectors in $V$. 
By theorem \ref{finiteG} the representations listed in the table have 
property $(\B_k)$. 
This together with lemma \ref{fix} and proposition \ref{directsum} treats 
finite reflection groups completely.  

\begin{center}
\begin{tabular}{|l||c|c|c|}
\hline
$\rho : G \to O(V)$ & $d$ & $k$ & $|G|$\\
\hline
\hline
$A_n$, $n \ge 1$ & $n+1$ & $n+1$ & $(n+1)!$\\
\hline
$B_n$, $n \ge 2$ & $2n$ & $2n$ & $2^n n!$\\
\hline
$D_n$, $n \ge 4$ & $2n-2$ & $2n$ & $ 2^{n-1} n!$\\
\hline
$I_2^{n}$, $n \ge 5$ & $n$ & $n$ & $2 n$\\
\hline
$G_2$ & $6$ & $6$ & $12$\\
\hline
$H_3$ & $10$ & $12$ & $120$\\
\hline
$H_4$ & $30$ & $120$ & $14400$\\
\hline
$F_4$ & $12$ & $24$ & $1152$\\
\hline
$E_6$ & $12$ & $27$ & $51840$\\
\hline
$E_7$ & $18$ & $56$ & $2903040$\\
\hline
$E_8$ & $30$ & $240$ & $696729600$\\
\hline
\end{tabular}
\end{center}

\section{Property $(\B_k)$ for finite rotation groups} \label{rotgroups}

Let us denote by $C_2^n$ the cyclic subgroup of $O(2)$ generated by the
counterclockwise rotation of $\mathbb{R}^2$ through the angle $2 \pi/n$. 
Here we have $d = |G| = n$, and for any non-zero vector $v \in \mathbb R^2$ 
its isotropy group $G_v$ is trivial. By corollary \ref{cor},
\emph{all finite rotation groups $C_2^n$ in the plane have property $(\B_n)$}.

\begin{remark*}
This together with the result for dihedral groups $I_2^n$ in section 
\ref{reflectiongroups}
gives a complete discussion of all finite subgroups of $O(2)$. 
\end{remark*}

Next we consider finite rotation groups in $3$-dimensional space.

Let $P$ be a $2$-dimensional linear subspace in $\mathbb{R}^3$. 
Any rotation $R$ in $O(P)$ can be extended to a rotation in $O(3)$, 
by setting $Rx = x$ for all $x \in P^{\bot}$ and using linearity. 
By extending each transformation in a cyclic subgroup $C_2^n$ of $O(P)$ in
this fashion, we obtain a cyclic subgroup of rotations in $O(3)$, which
will be denoted by $C_3^n$.

On the other hand, if $S$ is a reflection in $O(P)$, then $S$ may also be
extended to a rotation in $O(3)$, in fact to the rotation through the angle
$\pi$ having the reflection line of $S$ in $P$ as its axis of rotation:
define $Sx = -x$ for all $x \in P^{\bot}$ and extend by linearity. By
extending each transformation in a dihedral subgroup $I_2^n$ of $O(P)$ to a
rotation in $O(3)$, the resulting set of rotations is a subgroup of $O(3)$
isomorphic with $I_2^n$; it shall be denoted $I_3^n$.

If $T$, $W$, and $H$ denote the subgroups of rotations in $O(3)$ which
leave invariant the regular tetrahedron, cube, and icosahedron each with
center in the origin, then the following list provides a complete
characterization of finite rotation groups in $\mathbb{R}^3$ (e.g.
\cite{GroveBenson}):
\[C_3^n, n\ge 1; I_3^n, n \ge 2; T; W; H.\]

\emph{Rotation groups of type $C_3^n$ have property $(\B_n)$} by
construction, since the linear subspace $P^{\bot}$ is left pointwise
invariant under the $C_3^n$-action, and on $P$ it restricts to the
$C_2^n$-action; so lemma \ref{fix} and the result for $C_2^n$ give the
statement.

\emph{Rotation groups of type $I_3^n$ have property $(\B_n)$}:  
Note first that here $d = n$. Moreover, we have the decomposition 
$\mathbb R^3 = P \oplus P^{\bot}$ into irreducible subrepresentations. 
In order to make $|G|/|G_{v_1}|$ minimal for $0 \ne v_1 \in P$ we may choose
$v_1$ to lie on some reflection line on $I_2^n$ in $P$; then 
$|G|/|G_{v_1}| = 2n/2 = n$.
For $0 \ne v_2 \in P^{\bot}$ we find that $|G_{v_2}|$ is the number of 
rotations (including the identity) in $I_2^n$. So $|G|/|G_{v_2}| = 2n/n = 2$. 
Application of corollary \ref{cor} gives the assertion.

\emph{The rotational symmetry group $T$ of the regular tetrahedron has
property $(\B_6)$}: 
We have $d = 6$ (e.g. \cite{CumminsPatera}, \cite{PateraSharpWinternitz}). 
Further the action of $T$ on $\mathbb R^3$ is irreducible. 
The elements of $T$ consist of rotations through 
angles of $2 \pi/3$ and $4 \pi/3$ about each of four axes joining vertices 
of the tetrahedron with centers of opposite faces, rotations through the angle 
$\pi$ about each of the three axes joining the midpoints of opposite edges, 
and the identity. So $|T| = 12$.
The isotropy groups of non-zero vectors on axes joining vertices with centers 
of opposite faces have cardinality $3$, those of non-zero vectors on axes 
joining the midpoints of opposite edges have cardinality $2$, and all other 
isotropy groups of non-zero vectors are trivial. 
Hence application of corollary \ref{cor} gives the statement.

\emph{The rotational symmetry group $W$ of the cube has
property $(\B_9)$}: 
We have $d = 9$ (e.g. \cite{PateraSharpWinternitz}). 
The action of $W$ on $\mathbb R^3$ is irreducible. The elements of $W$ consist 
of rotations through angles of $\pi/2$, $\pi$, and $3 \pi/2$ about each of 
three axes joining the centers of opposite faces, rotations through angles of $2 \pi/3$ and $4 \pi/3$ about each of four axes joining extreme opposite vertices, rotations through the angle $\pi$ about each of six axes joining midpoints 
of diagonally opposite edges, and the identity. Thus $|W|=24$.
The isotropy groups of non-zero vectors on axes joining the centers of opposite 
faces have cardinality $4$, those of non-zero vectors on axes 
joining extreme opposite vertices have cardinality $3$, 
those of non-zero vectors on axes 
joining midpoints of diagonally opposite edges have cardinality $2$,
and all other isotropy groups of non-zero vectors are trivial. 
Apply corollary \ref{cor}.

\emph{The rotational symmetry group $H$ of the regular icosahedron has
property $(\B_{15})$}:
We have $d = 15$ (e.g. \cite{CumminsPatera}). 
The action of $H$ on $\mathbb R^3$ is irreducible. The elements of $H$ consist 
of rotations through angles of $2 \pi/5$, $4 \pi/5$, $6 \pi/5$, and $8 \pi/5$ 
about each of the six axes joining extreme opposite vertices, rotations 
through angles of $2 \pi/3$ and $4 \pi/3$ about each of ten axes joining centers 
of opposite faces, rotations through the angle $\pi$ about each of fifteen 
axes joining midpoints of opposite edges, and the identity. 
Therefore $|H|=60$.
The isotropy groups of non-zero vectors on axes joining extreme opposite 
vertices have cardinality $5$, those of non-zero vectors on axes 
joining centers of opposite faces have cardinality $3$, 
those of non-zero vectors on axes 
joining midpoints of opposite edges have cardinality $2$,
and all other isotropy groups of non-zero vectors are trivial. 
Apply corollary \ref{cor}.

\vspace{0.3cm}

The following table collects the results for finite rotation groups in two 
and three dimensions obtained in this section. The groups
in the first column of the table are meant to stay for their standard
representation, $d$ is the integer associated to representations in
\ref{dD}, and $k$ is as in corollary \ref{cor}. 

\begin{center}
\begin{tabular}{|l||c|c|c|}
\hline
$\rho : G \to O(V)$ & $d$ & $k$ & $|G|$\\
\hline
\hline
$C_2^n$, $n \ge 1$ & $n$ & $n$ & $n$\\
\hline
$C_3^n$, $n \ge 1$ & $n$ & $n$ & $n$\\
\hline
$I_3^{n}$, $n \ge 2$ & $n$ & $n$ & $2 n$\\
\hline
$T$ & $6$ & $6$ & $12$\\
\hline
$W$ & $9$ & $9$ & $24$\\
\hline
$H$ & $15$ & $15$ & $60$\\
\hline
\end{tabular}
\end{center}

\begin{remark*}
Observe that in this table we always have $d=k$, i.e., the respective representation has 
property $(\B_d)$.  
Since we need at least regularity $C^d$ for a curve 
in the orbit space to be liftable once differentiably (theorem 
\ref{recall}), we cannot expect to improve these results. 
Evidently this remark applies also for those representations in the table 
in section \ref{reflectiongroups} with $d=k$. 
\end{remark*}

\vspace{1cm}

\end{document}